\documentclass[10pt]{amsart}
\usepackage{amsfonts,amsthm,amssymb,euscript}
\usepackage[pdftex]{graphicx}
\usepackage{graphics}
\graphicspath{ {./images/} }
\usepackage[matrix,arrow]{xy}
\usepackage[active]{srcltx}
\usepackage{caption}
\usepackage[backend=biber,style=numeric,maxnames=50,maxcitenames=50,giveninits=true,doi=false,isbn=false,url=false,eprint=false]{biblatex}
\usepackage{subfiles}

\usepackage[colorlinks=true]{hyperref}

\newcommand{\C}{\mathbb{C}}
\newcommand{\R}{\mathbb{R}}
\newcommand{\Z}{\mathbb{Z}}

\newcommand{\N}{\mathbb{N}}

\theoremstyle{plain}
\newtheorem{Theorem}{Theorem}[section]
\newtheorem{thm}[Theorem]{Theorem}
\newtheorem{lem}[Theorem]{Lemma}

\newtheorem{cor}[Theorem]{Corollary}
\newtheorem{con}[Theorem]{Conjecture}

\theoremstyle{definition}
\newtheorem{defi}[Theorem]{Definition}
\newtheorem{rem}[Theorem]{Remark}

\title[Quadratic growth of geodesics on the two-sphere]{Quadratic growth of geodesics on the two-sphere}
\author{Bernhard Albach}
\address[Bernhard Albach]{RWTH Aachen}
\email[Bernhard Albach]{\texttt{albach@mathga.rwth-aachen.de}}

\begin{document}

\begin{abstract}
We prove that for any reversible Finsler metric on $S^2$, the number of prime closed geodesics grows at least quadratically with respect to length. This result significantly improves the previously known growth rates, which were established by Hingston. The main tools are an improvement on Franks' theorem about the number of periodic points of area-preserving annulus maps, and the theory of cylindrical contact homology in the complement of a link.

\end{abstract}

\maketitle

\tableofcontents

\section{Introduction}

\subsection*{Historical overview}
Let $(M,g)$ be a closed Riemannian surface. The task of counting closed geodesics has historically been a driving force in geometry and calculus of variations. If the surface has a nontrivial fundamental group, Hadamard~\cite{Had} showed that by minimizing the length of loops one obtains a closed geodesic in each free homotopy class. Since these classes are in one-to-one correspondence with conjugacy classes of the fundamental group, one not only obtains infinitely many closed geodesics, but it is even possible to make statements about their growth rate with respect to length. Defining $P^t(g)$ to be the number of geometrically distinct closed geodesics with length at most equal to $t$, one gets that: 
\begin{enumerate}
\item If the genus of $M$ is larger than $1$ $\Rightarrow$ $\liminf_{t\rightarrow \infty}\frac{\log(P^t(g))}{t}>0$.
\item If the genus of $M$ is equal to $1$ $\Rightarrow$ $\liminf_{t\rightarrow \infty}\frac{\log(P^t(g))}{\log(t)}\geq 2$.
\end{enumerate}

The case of $S^2$ is considerably more difficult, and many different techniques have been employed to tackle the closed geodesics problem, going back to Poincaré \cite{P1905} who, while studying the three-body problem, outlined a proof of the existence of at least one non-constant closed geodesic on a convex surface. The first general result was given by Lyusternik and Schnirelmann \cite{ballmann1978satz},\cite{Klingenberg2012},\cite{MR0029532},\cite{lyusternik1929probleme}, who showed the existence of at least three disjoint simple closed geodesics. The combined works of Bangert \cite{B1} and Franks \cite{MR1161099}
showed the existence of infinitely many closed geodesics. 

This result was further improved by Hingston \cite{NH}, who gave the first result for a lower bound on the growth rate of $P^t(g)$ by showing that it is at least the growth rate of the prime numbers: 
\begin{equation}
\label{old_growth_rate}
\liminf_{t \rightarrow \infty} P^t(g)\frac{\log(t)}{t}>0. 
\end{equation}
If the metric has a simple closed geodesic with inverse rotation number different from $1$, a result of Angenent \cite{MR2179729}, or alternatively a result of Hryniewicz-Momin-Salom\~ao \cite{UHAM}, shows that in this case the growth rate is much stronger: it is at least quadratic. 

More generally, it is proved in \cite{UHAM} that if a Reeb flow on the tight $S^3$ has a pair of periodic orbits forming a Hopf link whose transverse rotation numbers are not inverse to each other, then the growth rate of geometrically distinct periodic orbits with respect to their periods is also at least quadratic.  

\subsection*{Main result}

When $g$ is a Finsler metric, $P^t(g)$ is defined just as before. The main result of this paper is the following improvement on the growth obtained by Hingston.

\begin{thm}
\label{mainint}
Let $g$ be a reversible Finsler metric on $S^2$. Then the growth rate of $P^t(g)$ is at least quadratic:
\begin{equation}
\label{main_quad_growth}
\liminf_{t\rightarrow \infty}\frac{\log(P^t(g))}{\log(t)}\ge 2.
\end{equation}
\end{thm}

The main points to be made are that the growth rate~\eqref{main_quad_growth} is much stronger than the known growth rate~\eqref{old_growth_rate} due to Hingston, and, even more importantly, that it is expected to be the slowest growth that holds for all Finsler-reversible two-spheres. No genericity conditions on the metric are required. 

\subsection*{Method}
If the metric is Riemannian, the Lyusternik-Schnirelmann theorem states that there exist at least three closed, simple geodesics. By the work of Bangert \cite{B1}, we can differentiate the following two cases: Either there exists a closed simple (without self-intersections) geodesic whose Birkhoff annulus is a global surface of section, or there exist two simple closed geodesics that do not intersect each other. This was generalized to reversible Finsler metrics by the combined works of De Philippis-Marini-Mazzucchelli-Suhr \cite{MR4696545} and Contreras-Knieper-Mazzucchelli-Schulz \cite{CKMSS}.

In the first case, there exists a well-defined return map, and Theorem \ref{mainint} reduces to the following statement, which is an improvement on Franks' theorem about periodic points of area-preserving maps:

\begin{thm}
\label{FI}
Let $A = (0,1)\times S^1$ be the open annulus with its standard area-form, and let $f\colon A \rightarrow A$ be an area-preserving homeomorphism isotopic to the identity. Assume that $f$ has some periodic point. Define $P^t(f)$ to be the number of geometrically distinct periodic orbits of $f$ with period at most $t$.
Then
\[
\liminf_{t\rightarrow \infty}\frac{\log(P^t(f))}{\log(t)}\ge 2.
\]
\end{thm}


For the proof, the following two statements are used:

\begin{thm}[Franks~\cite{MR0951509},\cite{MR1161099},\cite{MR0967632}]
\label{franks_thm_twist}
Let $f \colon A \rightarrow A$ be an area-preserving diffeomorphism isotopic to the identity, $\pi\colon(0,1)\times \R\rightarrow A$ be the projection from the universal cover of $A$ and $F\colon (0,1)\times \R \rightarrow (0,1)\times \R$ be a lift of $f$ with respect to~$\pi$. Denote by $p_2:(0,1)\times\R\to\R$ the projection onto the $\R$-factor. If there exist $z_1,z_2 \in (0,1) \times \R$ and $p\in \Z$, $q\in \N$ coprime such that:
\begin{enumerate}
    \item $$ \liminf_{n\rightarrow \infty} \frac{p_2\circ F^n(z_1)}{n}\le \frac{p}{q}\le \limsup_{n\rightarrow \infty} \frac{p_2\circ F^n(z_2)}{n}, $$
    \item $\pi(z_1)$ and $\pi(z_2)$ are recurrent points of $f^m$ for every $m\in \Z$,
\end{enumerate}
then $f$ has a periodic point $z_0$ of period $q$ such that $$ \lim_{n\rightarrow \infty} \frac{p_2\circ F^n(z)}{n}=\frac{p}{q} $$ holds for any $z\in \pi^{-1}(z_0)$.
\end{thm}

\begin{lem}\label{numbertheorey}
Let $(a,b) \subset \R$ be a non-empty open interval and let $$ P^t(a,b) \:= \# \left\{(p,q) \ \colon \ p\in \Z, q\in \N, \frac{p}{q}\in (a,b), \text{$p$ and $q$ are coprime}, q\le t \right\}. $$ Then, 
\[\lim_{t\rightarrow \infty}\frac{\log(P^t(a,b))}{\log(t)} = 2\]
holds.
\end{lem}

In general, some maps do not satisfy the twisting condition of Franks' theorem, such as irrational rotations, but these also have no periodic points. To prove Theorem~\ref{FI}, we use the existence of a fixed point to show that the dynamical system must fulfill a twisting condition, albeit possibly for a different annulus. More precisely, we collapse the boundary to obtain, up to taking an iterate, a lift $\tilde{f}:S^2\rightarrow S^2$ of $f$ that has at least three fixed points. By a theorem of Franks, this map must have infinitely many periodic points. If there are infinitely many fixed points, the theorem is proved because the growth rate as defined above would be infinite. Thus, we will assume that there exists a periodic point with a period greater than one. We consider an isotopy from the identity to $f$ for which a result from Le Calvez \cite{MR2217051}, \cite{MR2059431} ensures the existence of a transverse foliation. For the precise definition, see Section \ref{annulusmap}. The trajectory under the isotopy of the periodic point now has to intersect some leaf of the transverse foliation. Removing the endpoints of this leaf gives a map $\hat f:A\rightarrow A$ that still has at least one fixed point. This fixed point and the periodic point give the desired twisting condition needed to apply Franks' Theorem~\ref{franks_thm_twist}. Quadratic growth of periodic points of $\hat{f}$ implies the same for $f$. Then Lemma~\ref{numbertheorey} completes the proof.

In the second case, there exist two disjoint simple closed geodesics. Using the Hilbert contact form on the unit tangent bundle of $S^2$, we lift the setting from a metric $g$ on $S^2$ to a contact form $\lambda$ on $S^3$. These two closed geodesics lift to a four-component link $L_1$ of closed Reeb orbits in $S^3$. In this case, the main theorem \ref{mainint} follows from combining the following result with Lemma \ref{numbertheorey}.

\begin{thm}
Let $\lambda$ be a tight contact form on $S^3$ that realizes a link $L$ that is transversely isotopic to $L_1$ as a link of closed Reeb orbits. Let $(a,b)\subset \R$ be any interval. Then there exists a $c\in \R_{>0}$, only depending on $(a,b)$, such that for each $(p,q)\in \Z\times \N$ coprime satisfying $\frac{p}{q}\in (a,b)$ there exists a closed Reeb orbit in homotopy class $y_{(p,q)}$ whose action is bounded from above by $cq$.
\end{thm}

A precise definition of the link $L_1$ and of the free homotopy classes $y_{(p,q)}$ of loops in the complement of $L_1$ will be given in Section \ref{homo}. The strategy of the proof in this case is to first choose a model contact form $\lambda_m$ and explicitly calculate its action-filtered cylindrical contact homology $HC_\ast^{y_{(p,q)},\le T}(\lambda_m,L)$ in the free homotopy classes $y_{(p,q)}$ of the complement of~$L$, as defined in~\cite{Momin}. Then, we use a neck stretching argument to show that $\lambda$ also has to have a periodic Reeb orbit in this homotopy class. This will also give us qualitative information about the existing geodesics. In this step, the form of the link will be essential to control limits of sequences of holomorphic curves. Our model contact form will be defined by the lift of the Hilbert contact form of a certain sphere of revolution. Using an explicit formula for a return map to a Birkhoff annulus of a special geodesic, we are able to determine all closed geodesics in the chosen homotopy classes. After a Morse-Bott perturbation, we obtain a family of contact forms on $S^3$ whose filtered homology we can calculate up to an arbitrarily large filtration value. 

\subsection*{Outlook}
In our theorem, we assumed the Finsler metric to be reversible, which is essential, as Katok \cite{Katok2},\cite{Ziller_1983} constructed examples of non-reversible Finsler metrics with only two closed geodesics. This example of Katok correlates with an important conjecture for $3$-manifolds: 

\begin{con}
Let $(M,\alpha)$ be a closed connected contact $3$-manifold. Then the associated Reeb flow has either exactly two or infinitely many periodic orbits.
\end{con}

This conjecture was recently established by Cristofaro-Gardiner, Hryniewicz, Hutchings and Liu for all contact manifolds whose associated contact structure has torsion first Chern class, see \cite{2or}. It is also known to be valid without any torsion assumption when the contact form is non-degenerate, due to a result by Colin, Dehornoy and Rechtman~\cite{MR4549092}. There are strong quantitative statements for the case where there are only two closed Reeb orbits, see Cristofaro-Gardiner, Hryniewicz, Hutchings and Liu \cite{CHHLC}. In light of the results of this paper, it would be interesting to know whether the following conjecture is true.

\begin{con}[Hryniewicz]
Let $(M,\alpha)$ be a closed connected contact $3$-manifold. Define $P^t(\alpha)$ to be the number of geometrically distinct periodic Reeb orbits of $\alpha$ with period at most $t$.
Then its Reeb flow has either exactly two periodic orbits or $$ \liminf_{t\rightarrow \infty}\frac{\log(P^t(\alpha))}{\log(t)}\ge 2. $$
\end{con}

Another question is 
to know whether the above lower bound on the growth rate is indeed sharp, as expected. 

\subsection*{Organization of the paper}
In Section 2, we present the necessary results regarding annulus maps. In Section 3, we recall the definition and some necessary facts about cylindrical contact homology in the complement of a link. Section 4 describes our model system, whose cylindrical homology we will compute in Section 5 after performing a Morse-Bott perturbation. Finally, in Section 6, we prove the main theorem. 

\subsection*{Acknowledgments}
I would like to thank my supervisor, Umberto Hryniewicz, for the helpful and enlightening discussions that made this work possible. I would also like to thank Patrice Le Calvez for his help in proving Theorem \ref{FI}. B.A. was partially funded by the Deutsche Forschungsgemeinschaft (DFG, German Research Foundation) – 320021702/GRK2326 – Energy, Entropy, and Dissipative Dynamics (EDDy). I am grateful to Michael Westdickenberg for the support during my PhD.

\section{Dynamics of annulus maps}\label{annulusmap}

Consider the open annulus $A = (0,1) \times S^1$. The goal of this chapter is to prove Theorem~\ref{FI}. This theorem is a consequence of the following statement.



\begin{thm}\label{anmagr2}
Let $f\colon S^2 \rightarrow S^2$ be an area-preserving homeomorphism isotopic to the identity. Then $f$ has either exactly two periodic points or 
\[
\liminf_{t\rightarrow \infty}\frac{\log(P^t(f))}{\log(t)}\ge 2.
\]
\end{thm}

In the above statement, if $f$ has exactly two periodic points, then these must be fixed points. In order to prove Theorem~\ref{anmagr2}, we need the following result from Franks~\cite{MR1161099}.

\begin{thm}[Franks~\cite{MR1161099}] \label{franks}
Let $f \colon A \rightarrow A$ be an area-preserving diffeomorphism isotopic to the identity, $\pi\colon(0,1)\times \R\rightarrow A$ be the projection from the universal cover of $A$ and $F\colon (0,1)\times \R \rightarrow (0,1)\times \R$ be a lift of $f$ with respect to $\pi$.

If there exist $z_1,z_2 \in (0,1) \times \R$ and $p\in \Z$, $q\in \N$ coprime such that:
\begin{enumerate}
    \item $$ \liminf_{n\rightarrow \infty} \frac{p_2\circ F^n(z_1)}{n}\le \frac{p}{q}\le \limsup_{n\rightarrow \infty} \frac{p_2\circ F^n(z_2)}{n} $$
    \item $\pi(z_1)$ and $\pi(z_2)$ are recurrent points of $f^m$ for every $m\in \Z$.
\end{enumerate}
Then $f$ has a periodic point $z_0$ of period $q$ such that $$ \lim_{n\rightarrow \infty} \frac{p_2\circ F^n(z)}{n}=\frac{p}{q} $$ holds for any $z\in \pi^{-1}(z_0)$.
\end{thm}





The next lemma shows that if the conditions of Theorem \ref{franks} are satisfied with a strict inequality, we have the desired growth rate:

\begin{lem}\label{number}
Let $(a,b) \subset \R$ be a non-empty open interval and define 
\[ P^t(a,b):=\# \left\{(p,q)\colon p\in \Z,\ q\in \N,\ \frac{p}{q}\in (a,b), \ p,q \:\mathrm{coprime}, \ q\le t \right\}. \]
Then
\[\lim_{t\rightarrow \infty}\frac{\log(P^t(a,b))}{\log(t)} = 2\]
holds.
\end{lem}

\begin{proof}
We start with the special case of $(a,b)=(0,1)$. In this case, the value $P^t(0,1)$ is equal to the sum $\sum_{n\in \mathbb{N} , n\le t}\phi(n)$, where $\phi(n)$ denotes the number of integers less than $n$ that are coprime to $n$. The function $\phi$ is also called Euler's totient function. The growth rate of this sum is well studied, and it is known that
\[ \sum_{n\in \mathbb{N}, n\le t}\phi(n) = \frac{3}{\pi^2}t^2 + O(t \ \log(t)), \]
see~\cite{MR3531455}. This shows that $\lim_{t\rightarrow \infty}\frac{\log(P^t(0,1))}{\log(t)}= 2$. Due to the properties of the logarithm function, translating or scaling the interval does not change the value $\lim_{t\rightarrow \infty}\frac{\log(P^t(a,b))}{\log(t)}$. This proves the lemma.
\end{proof}

In order to use Franks' theorem, we will produce from $f$ an area-preserving annulus homeomorphism with two points that satisfy $$ \lim_{n\rightarrow \infty} \frac{p_2\circ F^n(z_1)}{n}< \lim_{n\rightarrow \infty} \frac{p_2\circ F^n(z_2)}{n}. $$
We use the theory of transverse foliations, which we will recall next. For a more detailed description, see, for example, \cite{MR3787834} by Le Calvez and Tal.

\begin{defi}
Let $f$ be a homeomorphism of an oriented surface $\Sigma$ without boundary. An identity isotopy $I$ is a continuous path $\{f_t\}_{t\in [0,1]}$ in the space of homeomorphisms of $\Sigma$ equipped with the compact-open topology, connecting $f$ and the identity map. \\
The fixed point set $\mathrm{fix}(I):=\cap_{t\in [0,1]} \mathrm{fix}(f_t)$ of $I$ is defined as the set of points which are fixed points of all maps $f_t$, $t\in [0,1]$. We will also call these points rest points of the isotopy $I$.\\
The trajectory of a point $p$ is defined as the path $I(p):=t\mapsto f_t(p)$, $t\in [0,1]$, and by $I(p)^n$ we denote the concatenation of the trajectories $I(p),I(f(p)),...,I(f(p)^{n-1})$.\\
A singular isotopy is an identity isotopy defined on an open subset of $\Sigma$, which will be called $\mathrm{dom}(I)$, such that the complement of $\mathrm{dom}(I)$ in $\Sigma$, denoted by $\mathrm{sing}(I)$, is a subset of the fixed point set $\mathrm{fix}(f)$. \\
A singular isotopy $I$ is called a maximal isotopy of $f$ if there exists no point $p$ in $\mathrm{dom}(I)\cap \mathrm{fix}(f)$ whose trajectory $I(p)$ is homotopic to zero in $\mathrm{dom}(I)$ as a closed loop.
\end{defi} 

\begin{rem}
An isotopy $I=\{f_t\}_{t\in [0,1]}$ where the $f_t$ are homeomorphisms of the whole surface $\Sigma$ can still be interpreted as a singular isotopy, since $\Sigma$ is an open subset of itself. In this case, the singular set $\mathrm{sing}(I)$ is empty.
\end{rem}

\begin{defi}
Let $\Sigma$ be an oriented surface. An oriented foliation $\mathcal{F}$ defined on an open subset of $\Sigma$ is called a singular oriented foliation of $\Sigma$. The open subset where $\mathcal{F}$ is defined will be denoted by $\mathrm{dom}(\mathcal{F})$, its complement by sing$(\mathcal{F})$.
\end{defi}

\begin{defi}
Let $I=\{f_t\}_{t\in [0,1]}$ be a singular isotopy. A singular oriented foliation $\mathcal{F}$ is called transverse to $I$ if
\begin{enumerate}
\item $\mathrm{dom}(I)=\mathrm{dom}(\mathcal{F})$,
\item $\forall z\in \mathrm{dom}(I)$ the trajectory $I(z)$ is homotopic in $\mathrm{dom}(I)$, relative to the endpoints, to a path positively transverse to the foliation and 
\item this transverse path is unique up to equivalence, where two such paths $\gamma$ and $\tilde{\gamma}$ are said to be equivalent if there exists a holonomic homotopy between them. A holonomic homotopy between $\gamma$ and $\tilde{\gamma}$ is a continuous map $H:[0,1]\times[0,1]\rightarrow \mathrm{dom}(\mathcal{F})$ such that 
\begin{enumerate}
\item $H(t,0)=\gamma (t), H(t,1)=\tilde{\gamma}(h(t))$ for an increasing homeomorphism $h:[0,1]\rightarrow [0,1]$ and
\item for all $t,s_1,s_2\in [0,1]$ the equation $\phi_{H(t,s_1)}=\phi_{H(t,s_2)}$ holds, where $\phi_{H(t,s)}$ denotes the leaf containing the point $H(t,s)$.
\end{enumerate}
\end{enumerate}
\end{defi}

The following theorem connects these definitions.

\begin{thm}[Le Calvez \cite{MR2059431},\cite{MR2217051} ]\label{tra}
Let $I$ be a maximal isotopy of $f\colon \Sigma \rightarrow \Sigma$. Then a foliation $\mathcal{F}$ transverse to $I$ exists.
\end{thm}

\begin{proof}[Proof of theorem \ref{anmagr2}]
Assume that $f$ has at least three periodic points $z_0,z_1,z_2$. Taking iterates of $f$, we can assume that these points are fixed points. The map $f$ is isotopic to the identity and it is possible to choose an isotopy $I=\{f_t\}_{t\in [0,1]}$ between the identity and $f$ such that the fixed points $z_0,z_1,z_2$ are rest points of the isotopy (see, for example, \cite{MR0212840}). Since $f$ is area-preserving, a result by Franks \cite{MR1161099} states that $f$ has either two or infinitely many periodic points. Consequently, since $f$ has at least three fixed points, infinitely many periodic points must exist. If there exist infinitely many fixed points, the proof is obviously finished as the growth rate would be infinite; therefore, assume that there exist only finitely many fixed points. Consequently, there exists a periodic point $p$ with period $n>1$. We will now show how to get a maximal isotopy $\tilde{I}$ for $f$ starting with the isotopy $I$.

We achieve this through an iterative procedure. Define $I^0:=I$. If there exists a $z\in \mathrm{fix}(f)\setminus \mathrm{fix}(I)$ such that $I^0(z)$ is contractible in $S^2\setminus \mathrm{fix}(I)$, we perturb $I^0$ on a compact subset of $S^2\setminus \mathrm{fix}(I)$ to an identity isotopy $I^1$ of $f$ such that $z\in \mathrm{fix}(I^1)$. Then we keep repeating this process. Since $f$ has only finitely many fixed points, the process terminates after finitely many iterations. The resulting identity isotopy $\hat{I}$ of $f$, defined on the whole of $S^2$, has the property that there exists no point in $\mathrm{fix}(f)\setminus \mathrm{fix}(\hat{I})$ whose trajectory is contractible in $S^2\setminus \mathrm{fix}(\hat{I})$. Then the singular isotopy $\tilde{I}:=\hat{I}|_{S^2\setminus \mathrm{fix}(\hat{I})}$ is maximal. Note that $\mathrm{fix}(I) \subset \mathrm{fix}(\hat I) = \mathrm{sing}(\tilde I)$.

By Theorem \ref{tra}, there exists a singular foliation $\mathcal{F}$ transverse to this isotopy $\tilde{I}$. By the properties of the transverse foliation and the fact that the point $p$ is periodic with period $n>1$, $I(p)^n$ forms a loop which, relative to the endpoints in $\mathrm{dom}(\mathcal{F})$, is homotopic to a loop $\gamma$ that is positively transverse to the foliation. 

Choose a leaf of $\mathcal{F}$ that intersects $\gamma$. Since $\gamma$ is positively transverse to the foliation, this leaf must be homeomorphic to a line with two disjoint endpoints, see, for example, \cite{MR3787834}. These endpoints $e_1,e_2$ must be points in the singular set of the foliation, thus points in the singular set of $\tilde{I}$ and therefore rest points of $\hat{I}$.

Restrict $f$ to $S^2\setminus \{e_1,e_2\}$ to obtain a homeomorphism $h\colon (0,1)\times S^1\rightarrow (0,1)\times S^1$. The isotopy $\hat{I}|_{S^2\setminus \{e_1,e_2\}}$ is an isotopy between the identity and $h$.


Since $I$ has at least three rest points, $\mathrm{fix}(\hat{I})$ has at least three elements because of $\mathrm{fix}(I) \subset \mathrm{fix}(\hat I)$. Thus, $\hat{I}|_{S^2\setminus \{e_1,e_2\}}$ has at least one rest point, which we denote by $\tilde{z}_1$. Denote the lift of $h$ to $(0,1)\times \R$ induced by the isotopy $\hat{I}|_{S^2\setminus \{e_1,e_2\}}$ by $\tilde h\colon(0,1)\times \R \rightarrow (0,1) \times \R$. Then we have $\lim_{n\rightarrow \infty} \frac{p_2\circ \tilde h^n(\tilde{z}_1)}{n}=0$. Since the trajectory of $p$ algebraically intersects the chosen leaf a positive number of times and $p$ is a periodic point, there exists a $c\neq 0$ such that $\lim_{n\rightarrow \infty} \frac{p_2\circ \tilde h^n(p)}{n}=c$ holds, allowing us to use Franks' theorem \ref{franks}. Finally, the combination of Theorems \ref{franks} and \ref{number} completes the proof.  
\end{proof}

\section{Cylindrical contact homology}\label{ch}

Let $M$ be a closed $3$-manifold. Following Momin \cite{Momin}, we want to define 
\[ HC_\ast^{y,\le T}(\lambda, L), \]
the filtered cylindrical contact homology of the complement of a link in $M$, where
\begin{itemize}
\item $\lambda$ is a contact form on $M$,
\item $L$ is a link of periodic Reeb orbits of $\lambda$,
\item $y$ is a free homotopy class of loops in $M\setminus L$ and
\item $T\ \in (0,\infty]$.
\end{itemize}
This homology is well defined only under certain assumptions on $(\lambda, L, y, T)$.

$I)$ Topological assumptions:
\begin{enumerate}
\item $c_1(\xi)$ vanishes on toroidal classes, where $\xi:= \mathrm{Ker}(\lambda)$.
\item $y$ is primitive, i.e. there exists no loop $c$ in $M\setminus L$ and $k\ge 2$ such that $c^k$ represents the class $y$.
\item $\forall$ connected component $ K \subset L \: \mathrm{and} \: \forall n\in \N $: $ y_K\neq [K]^n $, where $y_K$ and $[K]$ are the free homotopy classes of loops in $(M\setminus L)\cup K$ induced by $y$ and $K$, respectively.
\item  $\forall$ connected component $ K \subset L \: \mathrm{and} \: \forall n\in \N$: $[K]^n\neq 0$.
\end{enumerate}

$II)$ Dynamical assumptions:
\begin{enumerate}
\item No closed Reeb orbit $\gamma$ in $M\setminus L$ with action $\le T$ is contractible in $M\setminus L$.
\item Every closed Reeb orbit $\gamma$ in $M\setminus L$ with action $\le T$ representing $y$ is non-degenerate. 
\end{enumerate}

By a closed, or periodic, Reeb orbit in $II)$ we mean an equivalence class of pairs $(x,T)$ where $T>0$ and $x:\mathbb{R}\to M$ is a $T$-periodic orbit of the Reeb flow of $\lambda$. We do not require $T>0$ to be the primitive period. Two pairs are equivalent if the underlying periodic orbits share a common point and have the same period. However, sometimes we may abuse notation and denote by the same symbol $P,\gamma\dots$ both the geometric image $x(\mathbb{R})$ or a parametrized loop $t\in \mathbb{R}/T\mathbb{Z} \mapsto x(Tt) \in M$ up to rigid rotations in the domain. By the free homotopy class of $\gamma$ we always mean that of a loop of the form $t \mapsto x(Tt)$. The set of periodic orbits of the Reeb flow of $\lambda$ in $M\setminus L$ representing $y$ and with action at most $T$ will be denoted by $P^{y,\leq T}(\lambda,L)$. The definition of a non-degenerate orbit will be given in~\ref{no-de}.

Unless otherwise stated, we assume for the remainder of this section that all tuples $(\lambda, L,y, T)$ satisfy these conditions. 

When $T=\infty$ one simply writes $HC_\ast^{y}(\lambda,L)$. In this case, it is possible to prove that $HC_\ast^{y}(\lambda,L)$ is independent of the choice of tuple $(\lambda, L,y,\infty)$ satisfying the above assumptions. When $T<\infty$ there are well-defined chain maps between corresponding filtered chain complexes, where a certain action shift needs to be taken into account.

\subsection{Conley-Zehnder index}
The chain groups of the homology will be generated by closed Reeb orbits. To obtain a graded homology, we require a grading for closed Reeb orbits, which will be given by the Conley-Zehnder index. We will now define this index:

\begin{defi}
\[ \Sigma^{\ast}:=\{\nu\colon [0,1]\rightarrow Sp(2) \ \text{piecewise smooth} \ \colon \ \nu(0)=I, \ \det \nu(1)-I \neq 0 \}. \]
Here $Sp(2)$ denotes the group of $2 \times 2$ real symplectic matrices.
\end{defi}

\begin{thm}
There exists a unique map $\mu \colon \Sigma^\ast \rightarrow \mathbb{Z}$ such that:
\begin{itemize}
\item Homotopy: The map is invariant under homotopy of paths in $\Sigma^{\ast}$.
\item Maslov index: If $\gamma\colon [0,1]\rightarrow Sp(2)$ is a loop with $\gamma(0)=\gamma(1)=I$, $\nu \in \Sigma^{\ast}$ and $M(\gamma)$ denotes the Maslov index of $\gamma$, then
\[\mu (\gamma \phi)= \mu(\phi)+ 2 M(\gamma).\]
\item Inverse: Given $\nu \in \Sigma^{\ast}$, define $\nu^{-1}(t):= \nu(t)^{-1}$. Then $\mu(\nu^{-1})=-\mu(\nu)$.
\item Normalization: $\mu(t\mapsto e^{\pi i t})=1$.
\end{itemize}
\end{thm} 

This map is called the Conley-Zehnder index. 

\begin{defi}\label{no-de}
A periodic Reeb orbit $\gamma$ of period $T$ is called non-degenerate if $1$ is not in the spectrum of $d\phi_T$ along $\gamma$, and a contact form $\lambda$ is called non-degenerate up to action $S$ if all closed orbits with action less than or equal to $S$ are non-degenerate.
\end{defi}

It can be used to grade closed Reeb orbits as follows. Let $\phi_s$ be the Reeb flow of the contact form $\lambda$ and let $\gamma$ be a periodic orbit with period $T$, not necessarily the primitive period. Then $d\phi_{Ts}|_\xi$ along $\gamma$, $s\in [0,1]$, is a family of symplectic maps which can be represented by a family of symplectic matrices $\varphi_s$ after choosing a $d\lambda$-symplectic trivialization $\beta$ of $\gamma^{\ast}\xi$. Consequently, non-degenerate periodic orbits together with such trivializations $\beta$ can be graded by the Conley-Zehnder index
\begin{equation}
\mu_{CZ}(\gamma,\beta)=\mu(\varphi_s).
\end{equation}
Let $\beta$ be induced by a global symplectic trivialization of $(\xi,d\lambda)$, which we will assume from now on to be the case. Then, while the value $\mu_{CZ}(\gamma,\beta)$ does depend on the homotopy class of the global trivialization, by Condition $I.1)$ the difference of the values of two homotopic orbits does not depend on the homotopy class of the global trivialization. From now on we will suppress the dependence on the global trivialization and simply write $\mu_{CZ}(\gamma)$.

\subsection{Chain complex}
Using the Conley-Zehnder index, it is possible to define the graded chain complex:

\begin{defi}
Consider the set  $P^{y,\le T}(\lambda;L)$ consisting of periodic Reeb orbits~$P$ of $\lambda$ in $M\setminus L$ in the homotopy class $y$ with action $\int_P\lambda \le T$. The vector space over $\mathbb{Z}/2\mathbb{Z}$ freely generated by $P^{y,\le T}(\lambda;L)$ and graded by the Conley-Zehnder index will be denoted by
\[C^{y,\le T}_\ast(\lambda;L)=\underset{P\in P^{y,\le T}(\lambda;L)\colon \mu_{CZ}(P)=\ast}\bigoplus \mathbb{Z}/2\mathbb{Z}\cdot q_P\]
\end{defi}

\subsection{Differential}
Next, we want to define a boundary map.

\begin{defi}
Given an almost complex structure $J$ on a manifold $X$ and Riemannian surface $(\Sigma,i)$, a map $u\colon \dot{\Sigma}\rightarrow X$ is called $J$-holomorphic if 
\[ Tu\circ i = J\circ Tu . \]
\end{defi}

\begin{defi}
The symplectization of $(M,\lambda)$ is the symplectic manifold $(M\times \mathbb{R},d(e^t\lambda))$, where $t$ denotes the $\mathbb{R}$-coordinate. The projection from the symplectization to the manifold $M$ will be denoted by
\[ \tau\colon M\times \mathbb{R}\rightarrow M \qquad \qquad (x,t)\mapsto x. \]
On the symplectization there exists an $\mathbb{R}$-action $\{g_c\}_{c \in \mathbb{R}}$ given by
\[ g_c\colon \mathbb{R}\times (M \times\mathbb{R})\rightarrow M\times\mathbb{R}, \qquad \qquad g_c(x,t)=(x,t+c). \]
\end{defi}

A special class of $\mathbb{R}$-invariant almost complex structures on the symplectization can be constructed as follows:\\
The contact structure $\xi$ equipped with the symplectic form $d\lambda$ forms a symplectic vector bundle. The set of $d\lambda$-compatible complex structures on this bundle will be denoted by $\mathcal{J}(\xi)$. The dependence on $d\lambda$ is suppressed from the notation, but should not be forgotten. This set is non-empty, and contractible with the $C^\infty$-topology. We will denote the Reeb vector field of a contact form $\lambda$ by $X_\lambda$ and, given a $J\in \mathcal{J}(\xi,d\lambda)$, we define an $\mathbb{R}$-invariant almost complex structure $\tilde{J}$ on $M \times \mathbb{R}$ by
\begin{align*}
\tilde{J}(\partial_t)&=X_{\lambda}\\
\tilde{J}_{|\xi}&=J
\end{align*}
This almost complex structure is compatible with $d(e^t\lambda)$, and the set of almost complex structures $\tilde{J}$ arising in this way will be denoted by $\mathcal{J}(\lambda)$.

\begin{defi}
Consider a closed Riemann surface $(\Sigma,i)$ and denote by $(\dot\Sigma,i)$ the Riemann surface obtained by deleting finitely many points from $\Sigma$. A $\tilde J$-holomorphic curve $u:\dot\Sigma \to M\times \mathbb{R}$ in a symplectization is called a finite-energy curve if
\[ 0 < E(u)=\sup_{\phi\in\Lambda}\int_{\dot{\Sigma}}u^\ast d\lambda_\phi < \infty \]
where 
\[\Lambda :=\{\phi\colon \mathbb{R}\rightarrow [0,1] \mid \phi^\prime\ge 0\}\]
and $\lambda_\phi$ denotes the $1$-form defined by
\[ (a,x)\mapsto \phi(a)\lambda_{|x} . \]
The value $E(u)$ is referred to as the Hofer energy of $u$.
\end{defi}

In the following we identify $\mathbb{R}\times S^1 \simeq \mathbb{C} P^1 \setminus \{[0,1],[1,0]\}$ via the diffeomorphism $(s,t) \mapsto [e^{2\pi(s+it)},1]$. Use this diffeomorphism to pull the complex structure from $\mathbb{C} P^1$ back. Hence, we can see $\mathbb{R} \times S^1$ as the punctured Riemann sphere, with two punctures.

\begin{defi}
Let $\tilde{J}$ be an element of $\mathcal{J}(\lambda)$, and $P_0,P_1$ elements of $P^{y,\le T}(\lambda,L)$. The moduli space $M_{\tilde{J}}^{y,\le T}(P_0,P_1;L)$ is defined as the space of equivalence classes of finite-energy $\tilde{J}$-holomorphic curves $u=(u_M,u_\mathbb{R}): \mathbb{R}\times S^1\rightarrow M\times \mathbb{R}$ that
\begin{itemize}
\item[(i)] have one positive puncture asymptotic to $P_0$,
\item[(ii)] have one negative puncture asymptotic to $P_1$,
\item[(iii)] $u(\mathbb{R}\times S^1) \cap \tau^{-1}(L) = \emptyset$,
\end{itemize}
where two elements are said to be equivalent if they are holomorphic reparametrizations of each other.
\end{defi}

In $(i)$ and $(ii)$ it is meant that, with $P_0 = (x_0,T_0)$ and $P_1=(x_1,T_1)$, there are constants $t_0,t_1 \in S^1$ such that $u_M(s,t+t_0) \to x_0(T_0t)$, $u_\mathbb{R}(s,t) \to +\infty$  as $s\to +\infty$, and $u_M(s,t+t_1) \to x_1(T_1t)$, $u_\mathbb{R}(s,t) \to -\infty$ as $s\to -\infty$.

Moduli spaces as above are always considered with the quotient topology induced by the $C^\infty_{loc}$-topology on the space of maps $\mathbb{R}\times S^1\rightarrow M\times \mathbb{R}$. 


\begin{thm}
Consider the set $\mathcal{J}(\lambda)$ equipped with the $C^\infty_{loc}$-topology. There exists a Baire residual set $\mathcal{J}_{reg}(\lambda)\subset \mathcal{J}(\lambda)$ with the following property. Given~$\tilde{J}\in \mathcal{J}_{reg}(\lambda)$ and any $P_0,P_1 \in P^{y,\le T}(\lambda,L)$ the space $M_{\tilde{J}}^{y,\le T}(P_0,P_1;L)$ is a smooth manifold of dimension $\mu_{CZ}(P_0)-\mu_{CZ}(P_1)$.
\end{thm}


The above theorem includes the statement that if $\mu_{CZ}(P_0)-\mu_{CZ}(P_1)<0$ and $\tilde{J}\in \mathcal{J}_{reg}(\lambda)$ then $M_{\tilde{J}}^{y,\le T}(P_0,P_1;L)$ is the empty set. The validity of the above statement relies crucially on the fact that loops in class $y$ are automatically simple (not iterated), and this forces $\tilde{J}$-holomorphic maps representing cylinders in $M_{\tilde{J}}^{y,\le T}(P_0,P_1;L)$ to be automatically somewhere injective.

The $\mathbb{R}$-action $\{g_c\}$ on the symplectization defines an $\mathbb{R}$-action on the moduli space $M_{\tilde{J}}^{y,\le T}(P_0,P_1;L)$ which is free and smooth if $\mu_{CZ}(P_0)-\mu_{CZ}(P_1)>0$. The next main step in the construction of $HC_\ast^{y,\le T}(\lambda, L)$ is the following compactness result.

\begin{thm}
\label{thm_finite_count}
Given a $\tilde{J}\in \mathcal{J}_{reg}(\lambda)$ and $P_0,P_1\in  P^{y,\le T}(\lambda,L)$ satisfying the condition $\mu_{CZ}(P_0)-\mu_{CZ}(P_1)=1$, the space $M_{\tilde{J}}^{y,\le T}(P_0,P_1;L)/\mathbb{R}$ is finite.
\end{thm}

Details of the proof can be found in~\cite{Momin,UHAM}. Here we give a sketch. Consider a sequence of finite-energy cylinders $u_n$ representing a sequence in the $0$-dimensional manifold~$M_{\tilde{J}}^{y,\le T}(P_0,P_1;L)/\mathbb{R}$. First note that the SFT Compactness Theorem from~\cite{SFT} is available since by conditions $I)$ and $II)$ all cylinders $u_n$ stay at a positive distance from the link $L$. Hence, by the SFT Compactness Theorem, up to choice of a subsequence, it can be assumed that the sequence converges to a so-called holomorphic building. In the particular situation at hand, the limiting holomorphic building is a cylindrical building in the sense that each level consists of exactly one finite-energy cylinder with one positive and one negative puncture. The crucial point is that there are no bubbling-off of planes. Indeed, if a plane bubbles off, then there are two possibilities: either its asymptotic limit is a cover of a component $K$ of $L$, or its asymptotic limit is a cover of a periodic orbit in $M\setminus L$. In any case, the plane does not intersect $\tau^{-1}(L)$. Indeed, if such intersections existed, they would be isolated and count algebraically positively, thus forcing the cylinders $u_n$ also to intersect $\tau^{-1}(L)$. This is a contradiction. Having established that such a plane does not intersect $\tau^{-1}(L)$, we now argue to obtain a contradiction as follows. In the former case, we get a contradiction to condition $I.4)$. In the latter case, we get a contradiction to condition $II.1)$. Having established that the limiting holomorphic building is a cylindrical building, we claim that no level is asymptotic to a cover of a component $K$ of $L$: if that was the case then, by considering $n$ large enough, we would conclude that the $M$-component of the loop $t\mapsto u_n(s,t)$ would be close to a cover of $K$, thus implying that $y_K = [K]^m$ for some $m\in \mathbb{N}$. This contradiction to condition $I.3)$ establishes the claim. Hence, the limiting cylindrical building consists of cylinders in $(M\setminus L)\times \mathbb{R}$. Now use that $\tilde{J}\in \mathcal{J}_{reg}(\lambda)$ together with the fact that Fredholm indices behave additively to conclude that the limiting holomorphic building consists of exactly one cylinder representing an element of $M_{\tilde{J}}^{y,\le T}(P_0,P_1;L)/\mathbb{R}$ which is the limit of a subsequence of the original sequence. We proved that the moduli space is a compact $0$-dimensional manifold, hence finite.

Using the above theorem, the desired boundary map can be defined as follows: 

\[ \partial(\lambda,\tilde{J})_\ast\colon C^{y,\le T}_\ast(\lambda,L)\rightarrow C^{y,\le T}_{\ast-1}(\lambda,L) \]
\[ q_P\mapsto \underset{P^\prime\in P^{y,\le T}(\lambda,L) \colon \mu_{CZ}(P^\prime)=\ast-1}\sum \#_2(M_{\tilde{J}}^{y,\le T}(P,P^\prime;L)/\mathbb{R})q_{P^\prime}\]
where $\#_2$ denotes the count mod $2$.

\begin{thm}
\label{thm_d_square_is_zero}
If $\tilde{J}$ is an element of $\mathcal{J}_{reg}(\lambda)$, then $\partial(\lambda,\tilde{J})_{\ast-1}\circ \: \partial(\lambda,\tilde{J})_\ast=0$. 
\end{thm}

The argument shares some analogous steps to that used in the proof of Theorem~\ref{thm_finite_count}. As usual in Floer theory, a gluing-compactness scheme can be implemented. Two cylinders obtained from a term in the count for $\partial(\lambda,\tilde{J})_{\ast-1}\circ \: \partial(\lambda,\tilde{J})_\ast$ can be glued since these are somewhere injective and $J \in \mathcal{J}_{reg}(\lambda)$. One gets finite-energy cylinders in $(M\setminus L) \times \mathbb{R}$ connecting periodic Reeb orbits in class $y$, with Fredholm index~$2$. Such glued cylinders belong to a non-compact $1$-dimensional moduli space, and the initial ``broken trajectory'' identifies one of its ends. Arguing as in the sketch of the proof of Theorem~\ref{thm_finite_count}, one shows that the holomorphic building obtained by looking at the other end of this connected component of the moduli space must be a two-level cylindrical building consisting of cylinders again counted by the map $\partial(\lambda,\tilde{J})_{\ast-1}\circ \: \partial(\lambda,\tilde{J})_\ast$. All assumptions $I)$ and $II)$ are used for this, analogous to the proof of Theorem~\ref{thm_finite_count}. In particular, there is automatic somewhere injectivity of cylinders in $(M\setminus L) \times \mathbb{R}$ with asymptotic limits representing $y$. Hence, the mod $2$ count vanishes and Theorem~\ref{thm_d_square_is_zero} follows.

As mentioned before, the homology of the chain complex $(C^{y,\le T}_\ast(\lambda,L),\partial(\lambda,\tilde{J})_\ast)$ will be denoted by $HC^{y,\le T}_\ast(\lambda,L)$. 

\subsection{Chain map}\label{chain map}
The aim of this section is to define a chain map between the homologies $HC^{y,\le T}_\ast(\lambda_-, L)$ and $HC^{y,\le T}_\ast(\lambda_+, L)$ of two different contact forms $\lambda_-$ and $\lambda_+$ provided the following assumptions are satisfied:
\begin{enumerate}
\item $\xi=\mathrm{Ker} (\lambda_-) = \mathrm{Ker} (\lambda_+)$ as oriented bundles.
\item If $f\colon M\rightarrow \mathbb{R}$ is determined by $\lambda_+=f \lambda_-$, then $f(x)>1 \;\forall x\in M$.
\end{enumerate}
We write $\lambda_+ \succ \lambda_-$ to denote the second condition.

If $\lambda_+ \succ \lambda_-$ holds, there exists a function $h\colon M \times \mathbb{R} \rightarrow (0,\infty)$ such that there exists $R_-,R_+\in \mathbb{R}$ such that $R_+>R_-$ and
\[ h(x,t)=e^{t-R_+}f(x) \quad t\ge R_+, \qquad h(x,t)=e^{t-R_-} \quad t\le R_-, \qquad \partial_th>0  \ \text{everywhere}.\]
The $2$-form $\omega=d(h\lambda_-)$ is an exact symplectic form, and $(M\times \mathbb{R},\omega)$ will be called an exact symplectic cobordism between $\lambda_-$ and $\lambda_+$. It can be split into three parts:
\[W_+(\lambda_+)=M\times [R_+,\infty), \quad W(\lambda_-,\lambda_+)=M\times [-R_-,R_+] \quad \text{and} \quad W_-(\lambda_-)=M\times (-\infty,-R_-]. \]
 Given $\tilde{J}_+ \in \mathcal{J}(\lambda_+)$ and $\tilde{J}_- \in  \mathcal{J}(\lambda_-)$, the set $\mathcal{J}(\tilde{J_-},\tilde{J}_+)$ is defined as the set of almost complex structures $\overline{J}$ such that:
 \begin{enumerate}
     \item $\overline{J}$ coincides with $\tilde J_+$ in a neighborhood of $W_+(\lambda_+)$,
     \item $\overline{J}$ coincides with $\tilde J_-$ in a neighborhood of $W_-(\lambda_-)$,
     \item $\overline{J}$ is $\omega$-compatible.
 \end{enumerate}
The set of such almost contact structures for which $\tau^{-1}(L)$ is a complex submanifold will be denoted by $\mathcal{J}(\tilde{J_-},\tilde{J}_+:L)$.

\begin{defi}
In this setting, the energy of a $\overline{J}$-holomorphic curve is defined as \[E(u):= E_-(u)+E_0(u)+E_+(0),\] where the three terms are defined as follows:
\begin{align*}
E_-(u):=&\sup_{\phi\in\Lambda}\int_{u^{-1}(W_-(\lambda_-))}u^\ast d\lambda_{-_\phi}\\
E_+(u):=&\sup_{\phi\in\Lambda}\int_{u^{-1}(W_+(\lambda_+))}u^\ast d\lambda_{+_\phi}\\
E_0(u):=&\int_{u^{-1}(W(\lambda_-,\lambda_+))}u^\ast \omega\\
\end{align*}
\end{defi}

\begin{defi}
Let $\lambda_+,\lambda_-$ be two contact forms that satisfy $\lambda_+ \succ \lambda_-$, let $\tilde{J}_+\in \mathcal{J}(\lambda_+), \tilde{J}_-\in \mathcal{J}(\lambda_-)$, $\overline{J}\in \mathcal{J}(\tilde{J_-},\tilde{J}_+:L)$ be almost complex structures, and let $P_-\in P^{y,\le T}(\lambda_-,L),P_+\in P^{y,\le T}(\lambda_+,L)$ be closed Reeb orbits.

The moduli space $M_{\overline{J}}^{y,\le T}(P_+,P_-;L)$ is defined as the space of equivalence classes of finite energy $\overline{J}$-holomorphic curves $\mathbb{R}\times S^1\rightarrow M\times \mathbb{R}$ that
\begin{itemize}
\item have one positive puncture asymptotic to $P_+$,
\item have one negative puncture asymptotic to $P_-$,
\item and do not intersect $\tau^{-1}(L)$,
\end{itemize}
where again two elements are said to be equivalent if they are reparametrizations of each other. Furthermore, all cylinders in the moduli space $ M_{\overline{J}}^{y,\le T}(P_+, P_-; L)$ will automatically be somewhere injective, since $P_+$ and $P_-$ are elements of the homotopy class $y$. 
\end{defi}


\begin{thm}
Consider the set $\mathcal{J}(\tilde{J_-},\tilde{J}_+; L)$ equipped with the $C^\infty_{loc}$ topology. There exists a Baire residual set $\mathcal{J}_{reg}(\tilde{J_-},\tilde{J}_+; L)\subset \mathcal{J}(\tilde{J_-},\tilde{J}_+; L)$ such that for every $P_+\in P^{y,\le T}(\lambda_+,L)$ and $P_-\in P^{y,\le T}(\lambda_-,L)$ the space $M_{\overline{J}}^{y,\le T}(P_+,P_-;L)$ is a manifold of dimension $\mu_{CZ}(P_+)-\mu_{CZ}(P_-)$. 
\end{thm}

The above theorem includes the statement that if $\mu_{CZ}(P_+)-\mu_{CZ}(P_-)<0$ and $\overline{J}\in \mathcal{J}_{reg}(\tilde{J_-},\tilde{J}_+; L)$ then $M_{\overline{J}}^{y,\le T}(P_0,P_1;L)$ is the empty set. As before, the validity of the above theorem relies crucially on the fact that loops in class $y$ are automatically simple (not iterated), forcing $\overline{J}$-holomorphic maps representing cylinders in $M_{\overline{J}}^{y,\le T}(P_0,P_1;L)$ to be automatically somewhere injective.


\begin{thm}
Given a $\overline{J}\in \mathcal{J}_{reg}(\tilde{J_-},\tilde{J}_+:L)$, $P_+\in P^{y,\le T}(\lambda_+),P_-\in P^{y,\le T}(\lambda_-)$ with $\mu_{CZ}(P_+)-\mu_{CZ}(P_-)= 0$, the space $M_{\overline{J}}^{y,\le T}(P_+,P_-;L)$ is finite. 
\end{thm}

The proof is the same as the proof of Theorem~\ref{thm_finite_count}. The only important additional remark to the argument is that $L \times \mathbb{R}$ has a $\overline{J}$-invariant tangent space, so that intersection points between $L\times\mathbb{R}$ and $\overline{J}$-holomorphic maps have positive intersection index, when isolated.

Therefore, the following map is well defined:

\begin{defi}\label{Phi}
\[\Phi(\overline{J})_\ast\colon C^{y,\le T}_\ast(\lambda_+,L)\rightarrow C^{y,\le T}_{\ast}(\lambda_-,L)\]
\[q_P\mapsto  \underset{P^\prime\in P^{y,\le T}(\lambda_-,L)\colon \mu_{CZ}(P^\prime)=\ast}\sum\#(M_{\tilde{J}}^{y,\le T}(P,P^\prime,L))q_{P^\prime}\]
\end{defi}


\begin{thm}
If $\overline{J} \in \mathcal{J}_{reg}(\tilde{J_-},\tilde{J}_+;L)$, then
\[\Phi(\overline{J})_{\ast-1}\circ \partial(\lambda_+,\tilde{J}_+)_\ast-\partial(\lambda_-,\tilde{J}_-)_\ast \circ \Phi(\overline{J})_\ast=0\] holds.
\end{thm}

The proof is the same as the proof of Theorem~\ref{thm_d_square_is_zero}, again noting that $L \times \mathbb{R}$ has a $\overline{J}$-invariant tangent space.

\subsection*{Chain homotopy}
The goal of this last subsection is to explain why the map $\Phi(\overline{J})_\ast$ from \ref{Phi} is independent of the choice of $\overline{J}\in \mathcal{J}_{reg}(\tilde{J_-},\tilde{J}_+;L)$ on the level of homology. Let $\tilde{J}_+\in \mathcal{J}_{reg}(\lambda_+)$, $\tilde{J}_-\in \mathcal{J}_{reg}(\lambda_-)$ and $\overline{J}_0,\overline{J}_1$ be two different choices in $\mathcal{J}_{reg}(\tilde{J_-},\tilde{J}_+;L)$. We will now see why are the maps $\Phi(\overline{J}_0)$ and $\Phi(\overline{J}_1)$ chain homotopic.

\begin{defi}
The set $\overline{\mathcal{J}}(\overline{J}_0,\overline{J}_1;L)$ is defined as the set of smooth homotopies $\{J_t\}_{t\in [0,1]}$ in $\mathcal{J}(\tilde{J_-},\tilde{J}_+;L)$ from $\overline{J}_0$ to $\overline{J}_1$. Such homotopy is, by definition, required to be given by a smooth map defined on $[0,1]\times M\times \mathbb{R}$. Given $P_+\in P^{y, \le T}(\lambda_+;L)$ and $P_-\in P^{y,\le T}(\lambda_-;L)$ define
\[  M_{\{J_t\}}^{y,\le T}(P_+,P_-;L) := \{ (t,[u]) \ \colon \ t\in[0,1], \ [u] \in M_{J_t}^{y,\le T}(P_+, P_-;L) \}. \]
\end{defi}

Once again, for every $t\in [0,1] $ every element in $ M_{J_t}^{y,\le T}(P_+,P_-;L)$ is somewhere injective since $P_+$ and $P_-$ are elements of the primitive homotopy class $y$.


\begin{thm}
Consider $\overline{\mathcal{J}}(\overline{J}_0,\overline{J}_1:L)$ equipped with the $C^\infty_{loc}$ topology. There exists a Baire residual set $\overline{\mathcal{J}}_{reg}(\overline{J}_0,\overline{J}_1:L) \subset \overline{\mathcal{J}}(\overline{J}_0,\overline{J}_1:L)$ such that for every $P_+\in P^{y, \le T}(\lambda_+,L),$ $ P_-\in P^{y,\le T}(\lambda_-,L)$ with~$\mu_{CZ}(P_+)-\mu_{CZ}(P_-)+1=0$, and for every $\{J_t\} \in \overline{\mathcal{J}}_{reg}(\overline{J}_0,\overline{J}_1:L)$, the moduli space $M_{\{J_t\}}^{y,\le T}(P_+,P_-;L)$ is a smooth $0$-dimensional manifold and, furthermore, the $t$-component of all elements in this space is different from $0$ and $1$.
\end{thm} 

The argument is in complete analogy to that used to prove Theorem~\ref{thm_finite_count}. It uses a compactness arguments hinging on the SFT Compactness Theorem and on assumptions $I)$ and $II)$.

Thus, the following map is well defined: 
\[T(\{J_t\})\colon C^{y,\le T}_\ast(\lambda_+,L)\rightarrow C^{y,\le T}_{\ast+1}(\lambda_-)\]
\[q_P\mapsto \underset{P^\prime\in P^{y,\le T}(\lambda_-,L)\colon CZ(P^\prime)=\ast+1}\sum\#(M_{\{J_t\}}^{y,\le T}(P,P^\prime;L))q_{P^\prime}\]

The desired result is now given by:

\begin{thm}\label{6}
If $\{J_t\} \in \overline{\mathcal{J}}_{reg}(\overline{J}_0,\overline{J}_1;L)$ then
\[ \Phi(\overline{J}_0)_\ast-\Phi(\overline{J}_1)_\ast = T(\{J_t\})_{\ast-1}\circ \partial(\lambda_+,J_+)_\ast-\partial(\lambda_-,J_-)_{\ast+1}\circ T(\{J_t\}) \]
holds.
\end{thm}

The proof of the above theorem is in analogogy to that of Theorem~\ref{thm_d_square_is_zero}. It uses the usual gluing-compactness scheme in Floer theory. The main points here are the following. Firstly, the homotopy class $y$ is primitive, so all cylinders asymptotic to periodic orbit orbits in $P^{y,\leq T}(\lambda_\pm,L)$ are somewhere injective. Secondly, since the $\tilde{J}_\pm$,~$\overline{J}_0$,~$\overline{J}_1$ and~$\{J_t\}$ are generic, if $P_\pm \in P^{y,\leq T}(\lambda_\pm,L)$ satisfy $\mu_{CZ}(P_+)=\mu_{CZ}(P_-)=k$ then $M_{\{J_t\}}^{y,\le T}(P_+,P_-;L)$ is a smooth $1$-dimensional manifold with boundary, its boundary points corresponding to those elements with $t$-component equal to $0$ or $1$, i.e. to cylinders counted by the map $\Phi(\overline{J}_0)_k,\Phi(\overline{J}_1)_k$. The ends of this moduli space are in 1-1 correspondence with the pairs of cylinders counted by the map $T(\{J_t\})_{k-1}\circ \partial(\lambda_+,J_+)_k - \partial(\lambda_-,J_-)_{k+1}\circ T(\{J_t\})_k$. The latter statement is proved by a compactness analysis of holomorphic buildings provided by the SFT Compactness Theorem as in the proof of Theorem~\ref{thm_d_square_is_zero}, where assumptions $I)$ and $II)$ are used, in combination with a gluing argument which can be performed since all relevant cylinders are Fredholm regular.  Theorem~\ref{6} follows from these facts.

\subsection{An explicit chain map}\label{comm}
Let $\tilde{J}_+\in \mathcal{J}_{reg}(\lambda_+)$ and $c\in (0,1)$ be such that $\lambda_+$ fulfills conditions $II)$ not only up to action $T$, but up to action $T/c$. Then, $\lambda_-:= c \lambda_+$ also satisfies these conditions up action~$T$, and we consider an almost complex structure $\tilde{J}_-\in \mathcal{J}(\lambda_-)$ uniquely determined by requiting that $\tilde{J}_-|_{\xi} = \tilde{J}_+|_\xi$.

The diffeomorphism 
\[ \phi \colon M\times \mathbb{R}\rightarrow M\times \mathbb{R}, \qquad (x,t)\mapsto (x,c^{-1}t) \]
has the property that $\phi^\ast \tilde{J}_+=\tilde{J}_-$. Thus, the $\tilde{J}_-$-holomorphic cylinders $u_-$ are exactly the maps $\phi^{-1}\circ u_+$, where $u_+$ is a $\tilde{J}_+$-holomorphic cylinder, and there is a 1-1 correspondence between moduli spaces of $\tilde{J}_-$-holomorphic cylinders and the moduli space of $\tilde{J}_+$-holomorphic cylinders. It can also be shown that $\tilde{J}_-$ is an element of $\mathcal{J}_{reg}(\lambda_-)$. Moreover, the following map induces an isomorphism at the level of homology:
\[j_\ast\colon C^{\le T/c,y}_\ast(\lambda_+,L)\rightarrow C^{\le T,y}_{\ast}(\lambda_-,L) \qquad q_{(x,S)} \mapsto q_{(x(c^{-1}\cdot),cS)} \]
Composing with the inclusion 
\[i_\ast\colon C^{\le T,y}_\ast(\lambda_+,L)\rightarrow C^{\le T/c,y}_{\ast}(\lambda_+,L)\]
gives a map:
\[j_\ast\circ i_\ast\colon C^{\le T,y}_\ast(\lambda_+,L)\rightarrow C^{\le T,y}_{\ast}(\lambda_-,L)\]
The following theorem is as in~\cite[ Lemma~4.7]{UHAM}.

\begin{thm}
\label{gleich}
For every $\overline{J}\in \mathcal{J}_{reg}(J_-,J_+;L)$, the maps $j_\ast\circ i_\ast$ and $\Phi(\overline{J})_\ast$ from Definition~\ref{Phi} are chain homotopic.
\end{thm}

\subsection{Splitting}\label{split}
Let $\lambda_-,\lambda,\lambda_+$ be three contact forms satisfying $\lambda_-\prec \lambda \prec \lambda_+$ and $\tilde{J}_\pm\in \mathcal{J} (\lambda_\pm),\tilde{J}\in \mathcal{J} (\lambda), J_1\in \mathcal{J}(\tilde{J}_-,\tilde{J};L),J_2\in \mathcal{J}(\tilde{J},\tilde{J_+};L)$ be almost contact structures. Using the $\mathbb{R}$-action $\{g_c\}$, the following $1$-parameter family of almost complex structures $J_R$ (also denoted $J_1\circ_R J_2$), is well-defined for $R\gg1$ by:

\[J_R=\begin{cases} 
(g_{-R})^\ast J_2 \: \mathrm{on} \: W_+(\lambda)\\
(g_R)^\ast J_1 \: \mathrm{on} \: W_-(\lambda)
\end{cases}\]


There exists a biholomorphic almost complex structure $J_R^\prime$ which is an element of $\mathcal{J}(\tilde{J}_-,\tilde{J}_+;L)$. One way to construct $J_R^\prime$ is as follows. Consider a map
\[\psi_R\colon (M\times \mathbb{R})\rightarrow (M\times \mathbb{R}) \]
\[\psi_R(x,t)\mapsto (x,\varphi_R),\]
where $\varphi_R\colon \mathbb{R}\rightarrow \mathbb{R}$ is a function such that 
\[\varphi_R(t)=\begin{cases}
a+R & a \le -R-\epsilon\\
a-R & a \geq R+\epsilon
\end{cases}, \qquad \varphi_R^\prime(t) >0 \quad \mathrm{everywhere},\]
for $\epsilon$ small. Then $J^\prime_R$ given by $(\psi_R)_\ast(J_R)$ satisfies the requirements.

In this setting, the energy of a $J_R$-holomorphic curve $u$ is defined as \[E(u)=E_{\lambda_-}(u)+E_{(\lambda_-,\lambda)}(u)+E_{\lambda}(u)+E_{(\lambda,\lambda_+)}(u)+E_{\lambda_+}(u),\] where the summands are given by 
\begin{align*}
E_{\lambda_-}(u)\colon=&\sup_{\phi\in\Lambda}\int_{u^{-1}(W_-(e^{-R}\lambda_-))}u^\ast d\lambda_{-_\phi}\\
E_{\lambda}(u)\colon=&\sup_{\phi\in\Lambda}\int_{u^{-1}(W(e^{-R}\lambda,e^R\lambda))}u^\ast d\lambda_{\phi}\\
E_{\lambda_+}(u)\colon=&\sup_{\phi\in\Lambda}\int_{u^{-1}(W_+(e^{R}\lambda_+))}u^\ast d\lambda_{+_\phi}\\
E_{(\lambda_-,\lambda)}(u)\colon=&\int_{u^{-1}(W(e^{-R}\lambda_-,e^{-R}\lambda))}u^\ast e^R\omega\\
E_{(\lambda,\lambda_+)}(u)\colon=&\int_{u^{-1}(W(e^{R}\lambda,e^{R}\lambda_+))}u^\ast e^{-R}\omega,\\
\end{align*}
where $\omega$ denotes sympelctic forms as in Section \ref{chain map}.

\section{The model system}
In this section, we will construct an explicit contact form on $S^3$, which, after a Morse-Bott perturbation, will be used in our neck stretching argument later.

\subsection{Sphere of revolution}
A sphere of revolution is a $2$-dimensional manifold $S$ in $\R^3$ diffeomorphic to $S^2$ and invariant under rotation around the $z$-axis. It has exactly two intersection points with the $z$-axis; we will denote the one with the higher $z$-value by $p_N$ and the one with the lower $z$-value by $p_S$. Each sphere of revolution is uniquely defined by its intersection with the $(x,z)$-plane, which forms a smooth closed curve. Let $M$ be the length of this curve. The intersection can then be parametrized counter-clockwise by arc-length as a smooth closed curve 
\[\sigma \colon \R/M\Z \rightarrow \R^2\]
\[s\mapsto (r(s),z(s))\]
such that 
$\sigma(0)$ corresponds to $p_S$.

A parametrization of the complement of the two points $p_S$ and $p_N$ is given by:
\begin{equation}
\begin{split}
& \Phi \colon (0,M/2)\times \R/2\pi\Z \rightarrow S\setminus\{p_S,p_N\}   \\
&(s,\theta)\mapsto \big(r(s)\cos(\theta),r(s)\sin(\theta),z(s)\big)
\end{split}
\end{equation}

Given an $s_0 \in (0,M/2)$, the arc-length reparametrization of the curve 
\[
\theta\mapsto (r(s_0)\cos(\theta),r(s_0)\sin(\theta),z(s_0))
\]
will be called the parallel $P_{s_0}$.
Arc-length parametrized curves satisfying $\dot{\theta}\equiv 0$ are called meridians.

Restricting the Euclidean metric of $\R^3$ to a sphere of revolution $S$ turns $S$ into a Riemannian surface where the pullback of the metric is given by 
\begin{equation}\label{Metric}
r(s)^2d\theta^2+ds^2.
\end{equation}
This shows that meridians are always geodesics and that a parallel $P_{s_0}$ is a geodesic exactly if $s_0$ is a critical point of the function $r$.

\subsection{The model system}
\begin{figure}
\centering{
\resizebox{90mm}{!}{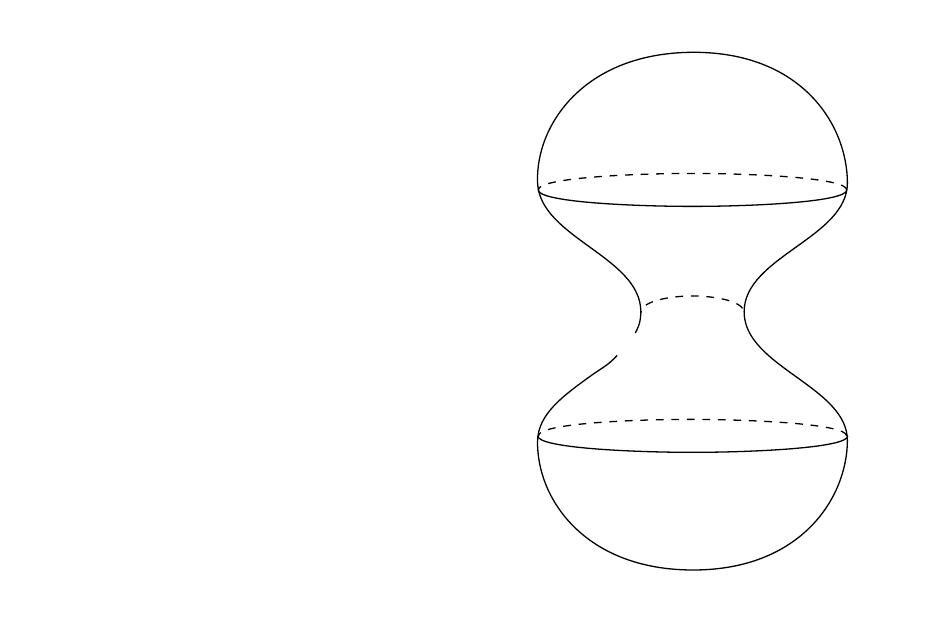}
\caption{The model sphere of rotation}
\label{Model}
}
\end{figure}

We define our model sphere of revolution $S_m$ by choosing a function $\sigma(s)=(r(s),z(s))$ such that
\begin{itemize}
\item $M> 2 \pi$,
\item $r(s)=\sin(s)$ for $s\in [0, \frac{\pi}{2}]$,
\item $\sigma$ is symmetric with respect to the $x$-axis and
\item the curve $r$ restricted to $(0,M/2)$ has three critical points; two maxima $s_{max_1},s_{max_2}$ at $s_{max_1}=\frac{\pi}{2}$ and $s_{max_2}=\frac{M}{2}-\frac{\pi}{2}$, which have the same value $r_{max}=1$ and one local minimum $s_{min}$, with value $r_{min}\in (0,1)$, at $s_{min}=M/4$.

\end{itemize}
See Figure \ref{Model}. The parallel $P_{s_{max_1}}$ will be called $D_+$, and the curve created by reversing the orientation will be denoted by $D_-$. Similarly, the parallel $P_{s_{max_2}}$ will be denoted by $E_+$, and the curve created by reversing the orientation will be denoted by $E_-$. The union $E_+\cup E_-\cup D_+ \cup D_-$ will be called $K_1$. Analogously, the parallel $P_{s_{min}}$ will be called $M_+$, and the curve created by reversing the orientation will be denoted by $M_-$.

\subsection{Unit tangent bundle}
The unit tangent bundle of 
$S_m$ with its Euclidean metric $g$ is defined as the $3$-dimensional manifold
\[ T^1S_m:=\{u\in TS_m: g(u,u)=1\}. \]
The projection $T^1S_m \rightarrow S_m$ will be denoted by $\pi$. An element $u \in T^1(S_m \setminus\{p_S,p_N\})$ is uniquely defined by its base-point $\Phi(s_0, \theta_0)\in S_m$ and the angle $\beta$ between the tangent vector of $P_{s_0}$ at $\Phi(s_0, \theta_0)$ and $u$.
Thus, a parametrization of $T^1(S_m\setminus\{p_S,p_N\})$ is given by:
\begin{equation}\label{cut}
\begin{split}
\tilde{\Phi}\colon (0,M/2)\times & \R/2\pi\Z\times \R/2\pi \Z \rightarrow T^1(S_m\setminus\{p_S,p_N\}) \\
&(s,\theta,\beta)\mapsto \cos(\beta)  \,\frac{\partial_\theta\Phi(s,\theta)}{r(s)} + \sin(\beta) \, \partial_s\Phi(s,\theta) 
\end{split}
\end{equation}

\begin{figure}
\centering{
\resizebox{30mm}{!}{
\begingroup%
  \makeatletter%
  \providecommand\color[2][]{%
    \errmessage{(Inkscape) Color is used for the text in Inkscape, but the package 'color.sty' is not loaded}%
    \renewcommand\color[2][]{}%
  }%
  \providecommand\transparent[1]{%
    \errmessage{(Inkscape) Transparency is used (non-zero) for the text in Inkscape, but the package 'transparent.sty' is not loaded}%
    \renewcommand\transparent[1]{}%
  }%
  \providecommand\rotatebox[2]{#2}%
  \newcommand*\fsize{\dimexpr\f@size pt\relax}%
  \newcommand*\lineheight[1]{\fontsize{\fsize}{#1\fsize}\selectfont}%
  \ifx\svgwidth\undefined%
    \setlength{\unitlength}{171.19931319bp}%
    \ifx\svgscale\undefined%
      \relax%
    \else%
      \setlength{\unitlength}{\unitlength * \real{\svgscale}}%
    \fi%
  \else%
    \setlength{\unitlength}{\svgwidth}%
  \fi%
  \global\let\svgwidth\undefined%
  \global\let\svgscale\undefined%
  \makeatother%
  \begin{picture}(1,0.929584)%
    \lineheight{1}%
    \setlength\tabcolsep{0pt}%
    \put(0,0){\includegraphics[width=\unitlength,page=1]{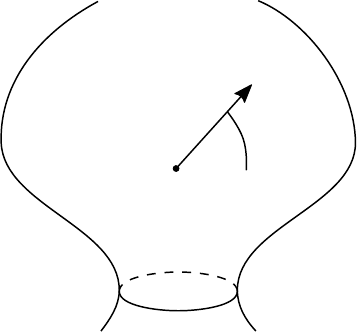}}%
    \put(0.47465594,0.51106408){\color[rgb]{0,0,0}\makebox(0,0)[lt]{\lineheight{1.25}\smash{\begin{tabular}[t]{l}$u$\end{tabular}}}}%
    \put(0,0){\includegraphics[width=\unitlength,page=2]{Winkel2.pdf}}%
    \put(0.5689339,0.48230067){\color[rgb]{0,0,0}\makebox(0,0)[lt]{\lineheight{1.25}\smash{\begin{tabular}[t]{l}$\beta$\end{tabular}}}}%
    \put(0.43902352,0.37879788){\color[rgb]{0,0,0}\makebox(0,0)[lt]{\lineheight{1.25}\smash{\begin{tabular}[t]{l}$\Phi(s_0, \theta_0)$\end{tabular}}}}%
  \end{picture}%
\endgroup%
}
\caption{The angle $\beta$}
\label{Winkel}
}
\end{figure}

\subsection{Hilbert contact form}

\begin{defi}
The contact form $\lambda_H$ on $T^1 S_m$ defined by 
\[(\lambda_{H})_u(v):=g(u,\pi_*v),u\in T^1S_m,v\in T_uT^1S_m\]
is called the Hilbert contact form and has the property that its Reeb flow coincides with the geodesic flow of $g$.
In the coordinates from parametrization \ref{cut} $\lambda_H$ is given by 
\begin{equation}\label{ctcform}
\lambda_H(s,\theta,\beta)=r(s)\cos(\beta)d\theta +\sin(\beta)ds
\end{equation}
\end{defi}

Its Reeb vector field is given by 
\begin{equation}\label{RVF}
X_{\lambda_H}(s,\theta,\beta)=\sin(\beta) \partial_s + \frac{r^\prime(s)}{r(s)}\cos (\beta) \partial_\beta +\frac{\cos(\beta)}{r(s)}\partial_\theta
\end{equation}
From this, we can see that the geodesic equations are:
\begin{equation}
\label{geoeq}
\begin{split}
s^\prime &= \sin(\beta)\\
\beta^\prime &= \frac{r^\prime(s)}{r(s)}\cos(\beta )\\
\theta ^\prime &= \frac{\cos(\beta)}{r(s)}
\end{split}
\end{equation}

\subsection{Clairaut integral}

\begin{defi}
The map
\begin{equation}
\begin{split}
K\colon T^1(S_m\setminus\{p_S,p_N\})\rightarrow \R\\
(s,\theta,\beta)\mapsto r(s)\cos(\beta)
\end{split}
\end{equation} 
is called the Clairaut integral.
\end{defi}

The important fact about it is that it is an integral of motion of the geodesic flow of the metric \eqref{Metric}, i.e. the value of the function is invariant under the geodesic flow. Thus, it is also invariant under the Reeb flow of $\lambda_H$. 

\begin{rem}
\label{fibr}
Let $S_m$ be the model form chosen previously and $k$ be a value in the interval $(-r_{min},r_{min})$. Then the level set $K^{-1}(k)$ is diffeomorphic to the torus $S^1\times S^1$. Thus, the set $K^{-1}(-r_{min},r_{min})$ is diffeomorphic to $(-r_{min},r_{min})\times S^1\times S^1$ and the Reeb vector field is tangent to the level sets $\{k\}\times S^1\times S^1$.
\end{rem}

\subsection{\texorpdfstring{Lift to  $\mathbf{S^3} $}{TEXT} }\label{lift}

\begin{figure}[h]
\centering{
\resizebox{80mm}{!}{
\begingroup%
  \makeatletter%
  \providecommand\color[2][]{%
    \errmessage{(Inkscape) Color is used for the text in Inkscape, but the package 'color.sty' is not loaded}%
    \renewcommand\color[2][]{}%
  }%
  \providecommand\transparent[1]{%
    \errmessage{(Inkscape) Transparency is used (non-zero) for the text in Inkscape, but the package 'transparent.sty' is not loaded}%
    \renewcommand\transparent[1]{}%
  }%
  \providecommand\rotatebox[2]{#2}%
  \newcommand*\fsize{\dimexpr\f@size pt\relax}%
  \newcommand*\lineheight[1]{\fontsize{\fsize}{#1\fsize}\selectfont}%
  \ifx\svgwidth\undefined%
    \setlength{\unitlength}{322.73727297bp}%
    \ifx\svgscale\undefined%
      \relax%
    \else%
      \setlength{\unitlength}{\unitlength * \real{\svgscale}}%
    \fi%
  \else%
    \setlength{\unitlength}{\svgwidth}%
  \fi%
  \global\let\svgwidth\undefined%
  \global\let\svgscale\undefined%
  \makeatother%
  \begin{picture}(1,0.93412466)%
    \lineheight{1}%
    \setlength\tabcolsep{0pt}%
    \put(0,0){\includegraphics[width=\unitlength,page=1]{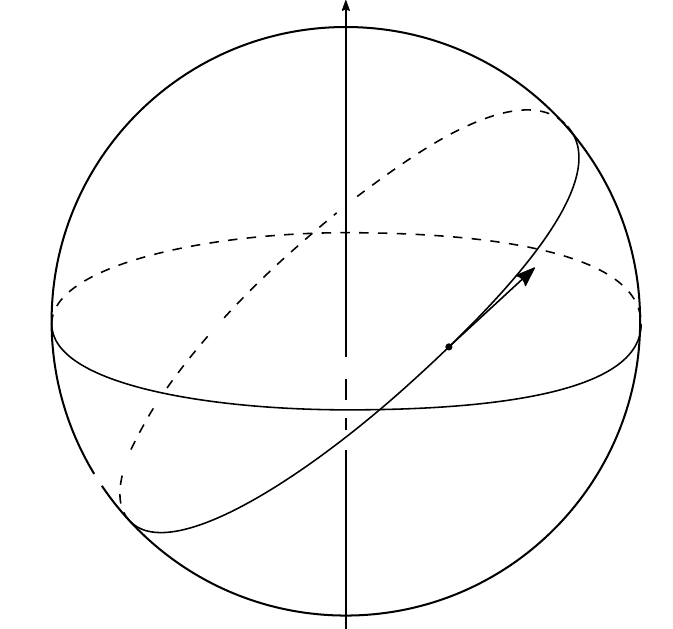}}%
    \put(0.68323023,0.40560878){\color[rgb]{0,0,0}\makebox(0,0)[lt]{\lineheight{0}\smash{\begin{tabular}[t]{l}$u$\end{tabular}}}}%
    \put(0,0){\includegraphics[width=\unitlength,page=2]{eulerangles.pdf}}%
    \put(0.55072864,0.29380538){\color[rgb]{0,0,0}\makebox(0,0)[lt]{\lineheight{1.25}\smash{\begin{tabular}[t]{l}$p$\end{tabular}}}}%
    \put(0,0){\includegraphics[width=\unitlength,page=3]{eulerangles.pdf}}%
    \put(0.43298775,0.46408226){\color[rgb]{0,0,0}\makebox(0,0)[lt]{\lineheight{1.25}\smash{\begin{tabular}[t]{l}$\phi (u)$\end{tabular}}}}%
    \put(0.53894902,0.40308559){\color[rgb]{0,0,0}\makebox(0,0)[lt]{\lineheight{1.25}\smash{\begin{tabular}[t]{l}$\theta (u)$\end{tabular}}}}%
    \put(0.60646663,0.33713647){\color[rgb]{0,0,0}\makebox(0,0)[lt]{\lineheight{1.25}\smash{\begin{tabular}[t]{l}$\nu(u)$\end{tabular}}}}%
    \put(0,0){\includegraphics[width=\unitlength,page=4]{eulerangles.pdf}}%
  \end{picture}%
\endgroup%
}
\caption{Euler angles}
\label{eulerangles}
}
\end{figure}

For every sphere of revolution $S$, there exists a double cover map from $S^3$ to $T^1S$. In order to explicitly define the map, we first introduce Euler angles:\\
We start with the round sphere $S^2$ of radius one. Let $u$ be an element of $T^1S^2\subset T\R^3$ which is not tangent to the equator.
In this case, $u$ defines an oriented great circle on $S^2$ that intersects the equator at exactly two different points. Let $p$ be the one where the tangent vector of the oriented great circle points into the upper hemisphere. We now associate three values $\phi(u),\theta(u)$ and $\nu(u)$ with $u$, see also Figure \ref{eulerangles}. \\
The value $\phi(u)$ is defined as the angle between the $x$-axis and $p$.\\
The value $\theta(u)$ is defined as the angle between $p$ and the basepoint of $u$ along the great circle defined by $u$.\\
The value $\nu(u)$ is defined as the angle between the tangent vector of the equator at $p$ and the tangent vector of the oriented great circle at $p$.

Since, conversely, these three values uniquely define a vector in $T^1S^2$ that is not tangent to the equator, the following map is a diffeomorphism:
\[\Phi_1\colon T^1S^2\setminus \{S^1 \sqcup S^1\}\rightarrow \R/2\pi\Z\times \R/2\pi\Z  \times (0,\pi)\]
\[u \mapsto (\phi(u),\theta(u),\nu(u))\]
Here $S^1 \sqcup S^1$ denotes the two-component link inside $T^1S^2$ given by the vectors tangent to the equator. Write $S^3$ as 
\[
S^3 = \{(r_1e^{it_1},r_2e^{it_2}),r_1,r_2\in [0,1],t_1,t_2\in \R/2\pi\Z:r_1^2+r_2^2=1\}\subset \C^2
\]
The map 
\[\Phi_2\colon S^3\setminus \{(z,0),(0,z):z\in S^1\} \rightarrow \R/2\pi\Z\times \R/2\pi\Z \times (0,\pi) \]
\[(r_1e^{it_1},r_2e^{it_2})\mapsto (t_1+t_2,t_1-t_2,2\arccos (r_1))\]
is a $2:1$ covering map, where the preimage of a given point $(\phi,\theta,\nu)$ consists of the two points $\pm(\cos(\frac{\nu}{2})e^{i\frac{\phi+\theta}{2}},\sin(\frac{\nu}{2})e^{i\frac{\phi-\theta}{2}})$. The composition $\Phi_1^{-1}\circ \Phi_2$ is a smooth $2:1$-covering map from $S^3\setminus \{(z,0),(0,z):z\in S^1\} $ to $T^1S^2\setminus \{S^1 \sqcup S^1\}$. This map can be smoothly extended to a $2:1$-covering map $l\colon  S^3 \rightarrow T^1S^2$ by mapping $(e^{it},0), t \in \R/2\pi\Z$ to twice the positively oriented equator of $S^2$ and $(0,e^{-it}),t \in \R/2\pi\Z$ to twice the negatively oriented equator.  

We want to use an orientation-preserving diffeomorphism $\hat{\Phi}\colon S_m \rightarrow S^2$ to extend the covering map $l$ to a covering map $\tilde{l}\colon S^3 \rightarrow T^1S_m$. Since $d\hat{\Phi}$ will generally not map $T^1S_m$ to $T^1S^2$, we have to use an additional rescaling, and we define the map 
\[\Phi_0\colon T^1S_m\rightarrow T^1S^2\]
\[\Phi_0(\eta)=\frac{d\hat{\Phi}(\eta)}{|d\hat{\Phi} (\eta)|_{g_{S^2}}},\]
where $g_{S^2}$ denotes the round metric on $S^2$. Then $\tilde{l}$ can be defined by $\tilde{l}=\Phi_0^{-1}\circ l$. 
We further require $\hat{\Phi}$ to map $M_+$ to the positively oriented equator of $S^2$. Using this covering map $\tilde{l}$, $(s,\theta,\beta)$ can be used as coordinates on $S^3$ (antipodal points have the same coordinates) with no fear of ambiguity.
The pullback of the Hilbert contact form $\lambda_H$ to $S^3$ will be denoted by $\lambda_m$ and in coordinates $(s,\theta,\beta)$ it still has the form $\lambda_m(s,\theta,\beta)=r(s)\cos(\beta)d\theta +\sin(\beta)ds$. We will denote the lifts to $S^3$ of the double covers of the velocity vectors of the curves $E_\pm,D_\pm$ and $M_\pm$ by $\tilde{E}_\pm,\tilde{D}_\pm$ and $\tilde{M}_\pm$, respectively.
Let $L_1 = \tilde{E}_\pm\cup\tilde{D}_\pm$. 
The reason why we lift a double cover of the curves is that otherwise the lifted curves would not be closed. 
The Clairaut integral also extends to a map $\tilde{K}\colon S^3\rightarrow \R$ by setting $\tilde{K}:= K\circ \tilde{l}$, and $\tilde{K}$ is an integral of motion of the Reeb flow of $\lambda_m$.

\begin{figure}
\centering{
\resizebox{90mm}{!}{
\begingroup%
  \makeatletter%
  \providecommand\color[2][]{%
    \errmessage{(Inkscape) Color is used for the text in Inkscape, but the package 'color.sty' is not loaded}%
    \renewcommand\color[2][]{}%
  }%
  \providecommand\transparent[1]{%
    \errmessage{(Inkscape) Transparency is used (non-zero) for the text in Inkscape, but the package 'transparent.sty' is not loaded}%
    \renewcommand\transparent[1]{}%
  }%
  \providecommand\rotatebox[2]{#2}%
  \newcommand*\fsize{\dimexpr\f@size pt\relax}%
  \newcommand*\lineheight[1]{\fontsize{\fsize}{#1\fsize}\selectfont}%
  \ifx\svgwidth\undefined%
    \setlength{\unitlength}{365.59364139bp}%
    \ifx\svgscale\undefined%
      \relax%
    \else%
      \setlength{\unitlength}{\unitlength * \real{\svgscale}}%
    \fi%
  \else%
    \setlength{\unitlength}{\svgwidth}%
  \fi%
  \global\let\svgwidth\undefined%
  \global\let\svgscale\undefined%
  \makeatother%
  \begin{picture}(1,0.56525801)%
    \lineheight{1}%
    \setlength\tabcolsep{0pt}%
    \put(0,0){\includegraphics[width=\unitlength,page=1]{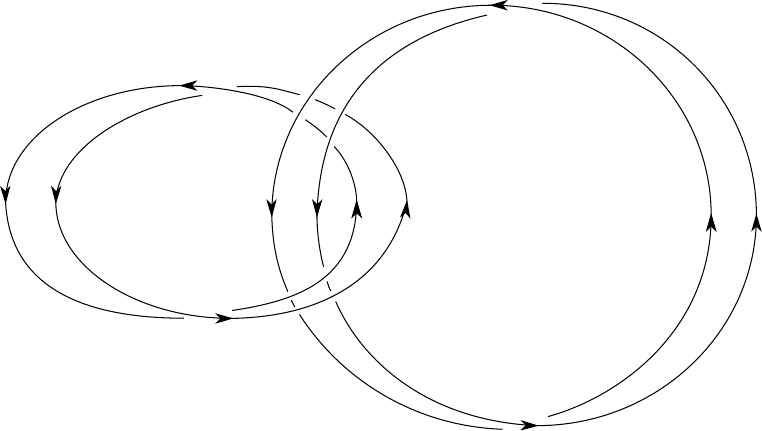}}%
    \put(0.02426027,0.42316396){\color[rgb]{0,0,0}\makebox(0,0)[lt]{\lineheight{1.25}\smash{\begin{tabular}[t]{l}$\tilde{E}_+$\end{tabular}}}}%
    \put(0.14697933,0.35280284){\color[rgb]{0,0,0}\makebox(0,0)[lt]{\lineheight{1.25}\smash{\begin{tabular}[t]{l}$\tilde{D}_+$\end{tabular}}}}%
    \put(0.53245694,0.46052757){\color[rgb]{0,0,0}\makebox(0,0)[lt]{\lineheight{1.25}\smash{\begin{tabular}[t]{l}$\tilde{D}_-$\end{tabular}}}}%
    \put(0.41477657,0.52231745){\color[rgb]{0,0,0}\makebox(0,0)[lt]{\lineheight{1.25}\smash{\begin{tabular}[t]{l}$\tilde{E}_-$\end{tabular}}}}%
  \end{picture}%
\endgroup%
}
\caption{The link $L_1$ projected to $\R^3$}
\label{Link}
}
\end{figure}


Using the Hurewicz homomorphism, we get $H_1(S^3\setminus L_1)\simeq \Z^4$. Furthermore, an explicit isomorphism is given by
\[ H_1(S^3\setminus L_1)\rightarrow \Z^4 \]
\[ [\gamma]\mapsto (lk(\gamma,\tilde{E}_+),lk(\gamma,\tilde{D}_+),lk(\gamma,\tilde{E}_-),lk(\gamma,\tilde{D}_-)), \]
where $\gamma$ is a continuous oriented curve whose homology class is denoted by $[\gamma]$.

\subsection{Satellites}

The notion of satellites will be crucial for us.


\begin{defi}
Let $S$ be an oriented surface, $\gamma \colon \R/2\pi\Z\rightarrow S$ be a regular curve, $\eta(t)$ a unit length normal vector along $\gamma$ such that $\{\dot{\gamma}(t),\eta(t)
\}$ is a positive basis of $T_{\gamma(t)}S$ for every $t$. Let $p\in \N,q \in \N_0$ be coprime numbers. A curve on $S$ will be called a $(p,q)-$satellite of $\gamma$ if it is homotopic, through immersed curves on $S$ that are nowhere tangent to $\gamma$, to the curve 
\[ \begin{split}
& \R/2\pi\Z\rightarrow S \\
t\mapsto \exp &_{\gamma(pt)}(\epsilon \sin(qt)\eta(pt)).
\end{split} \]

\end{defi}

\begin{lem}
\label{linking}
Let $\gamma_{(p,q)}$ be a $(p,q)$-satellite of the positively oriented equator $\gamma$ of $S^2$. The lift of $(\gamma_{(p,q)},\dot{\gamma}_{(p,q)})$ under the map $l$ is homotopic in $S^3\setminus l^{-1}(\gamma)$ to the curve
\[\R/2\pi \Z\rightarrow S^3\]
\[(c_1 e^{i\frac{p t}{2}},c_2e^{i\frac{(p-2 q)t}{2}}),\]
for $c_1,c_2\in \R_{>0}$ such that $c_1^2+c_2^2=1.$
\end{lem}

\begin{proof}
We argue similarly to \cite{MR2179729}. Define a curve $\tilde{\gamma}_{(p,q)}$ on the cylinder $$C:=\{(\alpha,z)\in \R^2\times \R \subset\R^3: |\alpha|=1\}$$ by $\tilde{\gamma}_{(p,q)}:\R/2\pi \Z \rightarrow C, t\mapsto (p \cdot t,\epsilon \sin (q\cdot t)),$ for an $\epsilon$ small. We denote the projection from the cylinder $C$ onto $S^2\setminus \{(0,0,1),(0,0,-1)\}\subset \R^3$ along the lines through the origin by $p\colon C\rightarrow S^2\setminus \{(0,0,1),(0,0,-1)\}$. The projection of $\tilde{\gamma}_{(p,q)}$ onto $S^2$ is then a $(p,q)-$satellite of the positively oriented equator. Since the Euler angles are defined using great circles on $S^2$, we study their behaviour under the projection $p$ now. The great circles are defined by intersections of planes containing the origin in $\R^3$ with $S^2$, and thus their preimage under the projection $p$ is defined by the intersection of the same plane with the cylinder $C$. Therefore, the preimage under the map $p$ of a great circle corresponding to the Euler angles $\phi$ and $\nu$ can be parametrized by the curve $\gamma_{(\phi,\nu)}\colon\theta  \mapsto (\theta+\phi,\sin(\nu)\sin(\theta))$. 
In order to get the Euler angles of the curve $t \mapsto (p\circ \tilde{\gamma}_{(p,q)}(t),\frac{d}{dt}(p\circ\tilde{\gamma}_{(p,q)}))$, we need to find values $\phi,\nu, \theta$, such that the tangent vector of $\gamma_{(\phi,\nu)}$ at the point $\theta$ agrees up to a scalar with the tangent vector $( \tilde{\gamma}_{(p,q)}(t),\frac{d}{dt}\tilde{\gamma}_{(p,q)})$. This is equivalent to solving 
\begin{align*}
(\tilde{\gamma}_{(p,q)}(t),\frac{d}{dt}\tilde{\gamma}_{(p,q)}(t))& =  (\gamma_{(\phi,\nu)}(\theta),(\frac{d}{dt}\gamma_{(\phi,\nu)})(\theta)) \\ 
\Leftrightarrow 
\begin{pmatrix}
    p \cdot t,\epsilon \sin (q\cdot t)\\
    p,\epsilon q \cos (q\cdot t)
\end{pmatrix}
&=
\begin{pmatrix}
    \theta +\phi,\sin(\nu)\sin(\theta)\\
    \theta'+\phi',\theta'\sin(\nu)\cos(\theta)+\nu'\cos(\nu)\sin(\theta)
\end{pmatrix}
\end{align*}
The solution is given by $\phi = (p-q)\cdot t,\theta =q\cdot t, \nu = \tilde{\epsilon},$ for $\tilde{\epsilon}$ small. Thus the Euler angles $\Phi_1(p\circ \tilde{\gamma}_{(p,q)},\frac{d}{dt}(p\circ\tilde{\gamma}_{(p,q)}))$ of the projection of $\tilde{\gamma}_{(p,q)}$ onto $S^2$ are given by the curve $t\mapsto ((p-q)\cdot t,q \cdot t,\tilde{\epsilon})$ in $\R/2\pi\Z\times \R/2\pi\Z  \times (0,\pi)$. The lemma then follows from the formula of the map $\Phi_2$ and the fact that every $(p,q)$-satellite of the equator $\gamma$ is homotopic to the curve $p\circ \gamma_{(p,q)}$.
 
\end{proof}
\begin{rem}
Given a satellite $\gamma_{(p,q)}$ of a curve on $S^2$ (or $S_m$), we can first lift this satellite to $T^1S^2$ (or $T^1S_m$) via the curve $(\gamma_{(p,q)},\dot{\gamma}_{(p,q)})$, and then to $S^3$ using the map $l$ (or $\tilde{l}$, respectively). We will call the resulting curve 'the lift of the satellite'. The previous lemma shows that, in general, the lifted curve need not be closed. In the case that it is not (which occurs if $p$ is odd), we instead lift a double cover of the satellite. By abuse of notation, we will still refer to the resulting curve as 'the lift of the satellite'.
\end{rem}

\subsection{Linking numbers}\label{linkingnumbers}

Using Lemma \ref{linking}, we can calculate the linking number of the lift $\tilde{\gamma}$ of the positively oriented equator $\gamma$ of $S^2$ and the lift of a satellite of $\gamma$. For this, we use the surface $\{(re^{it},\sqrt{1-r^2}):r\in [0,1],t\in [0,2\pi)\}$, which is a Seifert surface of $\tilde{\gamma}$. Using diffeomorphisms between $S_m$ and $S^2$, we
find that the linking numbers of $\tilde{E}_\pm$ and the lift $\tilde{\gamma}_{p,q}$ of a $(p,q)$-satellite $\gamma_{p,q}$ of $E_+$ are 
\begin{equation}\label{lkn}
\begin{split}
\begin{split}
lk(\tilde{\gamma}_{p,q},\tilde{E}_+)&=p-2q\\
lk(\tilde{\gamma}_{p,q},\tilde{E}_-)&=-p
\end{split}
\:\mathrm{if} \: p\: \mathrm{is\: odd},\:\\
\begin{split}
lk(\tilde{\gamma}_{p,q},\tilde{E}_+)&=\frac{p-2q}{2}\\
lk(\tilde{\gamma}_{p,q},\tilde{E}_-)&=-\frac{p}{2}
\end{split}
\:\mathrm{if} \: p\: \mathrm{is\: even}.\:
\end{split}
\end{equation} 
Using the same Seifert surface, we also get the following linking numbers:
\begin{equation}\label{linkingnumbersorbits}
\begin{split}
lk(\tilde{E}_+,\tilde{E}_-)=lk(\tilde{E}_+,\tilde{D}_-)&=lk(\tilde{D}_+,\tilde{E}_-)=lk(\tilde{D}_+,\tilde{D}_-)=-1\\
lk(\tilde{E}_+,\tilde{D}_+)&=lk(\tilde{E}_-,\tilde{D}_-)=+1
\end{split}
\end{equation}

Now that we know the linking numbers of satellites, we show that in our case all closed geodesics can be interpreted as satellites.
From equation \ref{geoeq} it follows that all geodesics $\gamma$ on $S$ can be divided into two cases:
\begin{itemize}
\item The derivative of the $\theta$ component is always zero.
\item The derivative of the $\theta$ component is never zero.
\end{itemize}
If $\theta^\prime$ is everywhere equal to zero, the curve is a meridian and it is easy to see that it is a $(1,1)-$satellite for $E_\pm$ and $D_\pm$. Otherwise, if $\theta^\prime$ is greater than zero, the closed geodesic is a satellite of $E_+$ or $D_+$, possibly of both, and if $\theta^\prime$ is less than zero, it is a satellite of $E_-$ or $D_-$, again possibly of both. The following theorem will give us more information about the type of satellite:

\begin{thm}[\cite{AABB}, Lemma 1.1]\label{wit}

Let $S_m$ be the model sphere of revolution, $\gamma$ be a geodesic on $S_m$ with Clairaut integral $K(\gamma) \notin \{0,\pm r_{min},\pm r_{max}\}$. Then there exist numbers $s_1,s_2\in (0,M/2)$ such that $r(s_1)=r(s_2)=|K(\gamma)|$, $\gamma$ is confined to the set \[\bigcup_{s\in [s_1,s_2]}P_s\]
and oscillates between $P_{s_1}$ and $P_{s_2}$, tangentially touching both alternately.
\end{thm}

Using this theorem, we can determine the type of satellite a closed geodesic is, depending on its Clairaut integral.\\
Let $\gamma$ be a closed geodesic with $K(\gamma)\in (0,r_{min})$. Then there are exactly two values $s_1$ and $s_2$ with $r(s_i)=K(\gamma)$ and by the previous theorem $\gamma$ oscillates between $P_{s_1}$ and $P_{s_2}$, tangentially touching both alternately. By the chosen form of the sphere of revolution $S$, it holds that $s_1<s_{max_1}$ and $s_2>s_{max_2}$. Thus, the curves $(D_+,\dot{D}_+)$ and $(E_+,\dot{E}_+)$ are homotopic to each other through immersed curves in $T^1S\setminus (\gamma,\dot{\gamma}) $. Therefore, $\gamma$ is the same type of satellite of both $D_+$ and $E_+$ and there exists a $(p,q)\in \N\times \N_0$ such that $\gamma$ is a $(p,q)$-satellite of $D_+$ and of $E_+$. 
Analogously, for $K\in (-r_{min},0)$ there exists $(p,q)\in \N\times \N_0$ such that $\gamma$ is a $(p,q)$-satellite of $D_-$ and of $E_-$.

In the case $K(\gamma) \in (r_{min},r_{max})$, the geodesic $\gamma$ is restricted to a strip around either $E_+$ or $D_+$. Therefore, there exists a $(p,q)\in \N\times \N_0$ such that $\gamma$ is a $(p,q)$-satellite of either $E_+$ or $D_+$. For the other one, it is a $(p,0)$-satellite.

Analogously, for $K(\gamma) \in (-r_{max},-r_{min})$ there exists a $(p,q)\in \N\times \N_0$ such that $\gamma$ is a $(p,q)$-satellite of either $E_-$ or $D_-$. For the other one, it is a $(p,0)$-satellite. 

If the closed geodesic $\gamma$ is a prime curve, then $p$ and $q$ are coprime. Using formula \ref{lkn} we can then compute the linking numbers of $\tilde{l}^{-1}(\gamma)$ and the curves $\tilde{E}_\pm$ and $\tilde{D}_\pm$. Since these linking numbers determine the homology class that $\tilde{l}^{-1}(\gamma)$ represents, we get the following lemma:

\begin{lem}\label{list}
Let $\gamma$ be a prime closed geodesic on the model sphere of revolution $S_m$ with $K(\gamma) \notin \{0,\pm r_{max}\}$. Then there exists $\tilde{p},\tilde{q} \in \N$ coprime such that $\tilde{\gamma}$, the lift of $\gamma$ to $S^3$, represents the following element in $H_1(S^3\setminus L_1)$:

\[
   \begin{array}{lr}
        (p-2q,p,-p,-p) \:\text{or} \:(p,p-2q,-p,-p) & \text{if } K(\gamma)\in (-r_{max},-r_{min})\\
        (-p,-p,p-2q,p-2q) & \text{if } K(\gamma)\in (-r_{min},0)\\
       (p-2q,p-2q,-p,-p), & \text{if } K(\gamma)\in (0,r_{min})\\
        (-p,-p,p-2q,p)\: \text{or}\: (-p,-p,p,p-2q) & \text{if } K(\gamma)\in (r_{min},r_{max})\\
        (\pm 1,\pm 1,\mp 1,\mp 1) & \text{if } K(\gamma)=\pm r_{min},
        \end{array}
  \]
where $(p,q)=(\tilde{p},\tilde{q})$ if $p$ is odd and $(p,q)=(\frac{\tilde{p}}{2},\frac{\tilde{q}}{2})$ if $p$ is even.

\end{lem}

\begin{rem}\label{rema}
Conversely, given a closed Reeb orbit in $S^3\setminus L_1$ representing an element $(a,b,c,d)$ in $H_1(S^3\setminus L_1)$, it is also possible to determine the range of the Clairaut integral of the closed geodesic it projects to. For example, assume $a=b,c=d$ and $(a,b,c,d)\neq (\pm 1,\pm 1,\mp 1,\mp 1)$. Then, according to the list above, the geodesic must have Clairaut integral in the range $(-r_{min},r_{min})$.
\end{rem}

\subsection{Existence of orbits}
Next, we will show the existence of certain closed geodesics using a surface of section.
Let $A$ be the Birkhoff annulus at $M_+$, which is defined as the set 
\[A:=\{(\Phi(s_{min},\theta),\beta)\in T^1S\colon \theta \in \R/2\pi \Z, \beta \in (0,\pi)\}.\] 
We will use the coordinate $\eta=- \cos(\beta)$ and parameterize the annulus by
\[\{(x,\eta), x \in \R/L\Z, \eta \in (-1,1)\},\]
where $L$ is the length of $P_{s_{min}}$.

The following theorem shows that the return map $\rho\colon A \rightarrow A$ is well-defined. 
\begin{thm}[\cite{AABB}, Lemma 2.1]\label{trans}
The forwards and backwards Reeb flow through any element of $A$ is transverse to $A$ and intersects $A$ again.
\end{thm}

Since the Clairaut integral stays constant under the Reeb flow, the map $\rho$ has to have the form $\rho(x,\eta)=(\tilde{\rho}(x,\eta),\eta)$ for a smooth function $\tilde{\rho}\colon A\rightarrow \R/L\Z$. Additionally, since $S$ is invariant under rotation around the $z$-axis, the map $\tilde{\rho}$ is also invariant under rotation around the $z$-axis and therefore the map $\rho$ can be written as
\begin{equation}\label{fmap}
     \rho(x,\eta)=(x+f(\eta),\eta),
\end{equation}

 for a smooth function $f\colon (-1,1)\rightarrow \R/L\Z$. We can choose the map $f$ such that $f(0)=0$, because the meridians are closed geodesics that intersect the Birkhoff annulus in the points $(x,0)$, $x\in \R/L\Z$. Then, by the symmetry of $S$, the function $f$ is an odd function.

\begin{defi}\label{windingnumber}
Given $u=(x,\eta)\in A$, the return time $\tau(u)$ is defined by \[\tau(u):= \min\{t\in \R_+: \phi_t(u)\in A \},\] where $\phi_t$ denotes the Reeb flow of $\lambda_H$. By Theorem \ref{trans} this value is always finite and because of the symmetry of $S$ the return time is independent of $x$. Therefore, we write it as $\tau(\eta)$. Also, by the symmetry of $S$, the return time $\tau(\eta)$ is an even function of $\eta$. \\
If $\eta \neq 0$, a geodesic intersecting the Birkhoff annulus in a point $(x,\eta)$ does not intersect the $z$-axis, and it is possible to define a winding number around the $z$-axis. For this purpose, we write the geodesic starting at the point $(x,\eta)\in A$ as 
\[t\mapsto (r(s(t))\cos(\theta(t)),r(s(t))\sin(\theta(t)),z(s(t))),t\in [0,\tau(\eta)].\]
The winding number is then defined by
\[W(u):=\frac{\theta(\tau(u))-\theta(0)}{2 \pi}.\]
The winding number is also independent of $x$ by the symmetry of $S$ and we write it as $W(\eta)$.\\
Given a closed geodesic $\gamma$ with minimal period $T$, intersecting the Birkhoff annulus $A$ in a point $(x,\eta)$ such that $\eta\neq 0$, we can also define its winding number $W(\gamma)$ by setting 
\[W(\gamma):=\frac{\theta(T)-\theta(0)}{2 \pi}.\] We will refer to this winding number as the winding number of the whole closed geodesic. It holds that $W(\gamma)=q\times W(\eta)$, where $q$ is the number of disjoint intersection points of the closed geodesic $\gamma$ and the Birkhoff annulus $A$.
\end{defi}

\begin{thm}[\cite{AABB} Lemmas 3.1 and 4.1]\label{tes}
Let $F$ be a primitive of $f$, i.e. $F^\prime=f$, such that $F(0)=M$, where $M$ is the length of the meridians.
Then, for all $\eta \in (-1,0)$ following equations hold:
\[f(\eta)=L W(\eta)-L \]
\[\tau (\eta)=F(\eta)-\eta f(\eta).\]

\end{thm}

Using this, we can prove the following theorem.

\begin{thm}\label{dec}
The map $f(\eta)$ from equation \ref{fmap} is strictly decreasing.
\end{thm}
\begin{proof}

We first prove this for $\eta \in (-1,0)$. Because $f$ is an odd function, the theorem then holds for all $\eta$.\\
Choose an $\eta$ in $(-1,0)$, $x\in \R/L\Z$ and write the geodesic $\gamma$ starting at $(x,\eta)$ as  
\[(r(s(t))\cos(\theta(t)),r(s(t))\sin(\theta(t)),z(s(t))).\]

This geodesic has a lot of symmetry, namely, it holds that
\[s(0)=s\left(\frac{1}{2}\tau(\eta)\right)=s(\tau(\eta))=M/4,\]
\[ s\left(\frac{1}{4} \tau(\eta)\right)=\max(s):=\max\{s(t),t\in[0,\tau(\eta)]\},\]
\[s\left(\frac{3}{4} \tau(\eta)\right)=\min(s):=\min\{s(t),t\in[0,\tau(\eta)]\}\]
and
\begin{equation}\label{rever}
s(t)=\begin{cases}
s\left(\frac{\tau(\eta)}{2}-t\right) & t \in \left[0,\frac{\tau(\eta)}{2}\right] \\
s\left(\tau(\eta)-\left(t-\frac{\tau(\eta)}{2}\right)\right) & t \in \left[\frac{\tau(\eta)}{2},\tau(\eta)\right]
\end{cases}\end{equation} 
These equations hold because of Theorem \ref{wit}, the fact that spheres of revolution are always invariant under reflections on planes containing the $z$-axis, and because we have chosen $S_m$ to be invariant under reflection on the $(x,y)$-plane.

Furthermore, by Theorem \ref{wit} and the geodesic equations \ref{geoeq}, the function $s$ is strictly decreasing on the interval $[\frac{1}{4} \tau(\eta),\frac{3}{4} \tau(\eta)]$. Thus, there exists a well-defined inverse function \[s^{-1}\colon [\min(s),\max(s)] \rightarrow \left[\frac{1}{4} \tau(\eta),\frac{3}{4} \tau(\eta)\right].\]  We can now compute the winding number:
\begin{align*}
2 \pi W(\eta)=\theta(\tau(\eta))-\theta(0)&=\int_0^{\tau(\eta)}\frac{\partial \theta}{\partial t}(t) dt\stackrel{(\ast)}{=}2 \int_{\frac{1}{4} \tau(\eta)}^{\frac{3}{4} \tau(\eta)}\frac{\partial \theta}{\partial t}(t) dt\\
&=2\int^{\min(s)}_{\max(s)}\frac{\partial \theta}{\partial t}(s^{-1}(z))\frac{\partial s^{-1}}{\partial z}(z)dz
\end{align*}

The third equality $(\ast)$ holds because of Equation \ref{rever}. Using the geodesic equation \ref{geoeq}, we get
\[ \pi W(\eta)=\int^{\min(s)}_{\max(s)}\frac{\cos(\beta(s^{-1}(z)))}{r(z)}\frac{1}{\sin(\beta(s^{-1}(z)))} dz\]
By the invariance of the Clairaut integral we get the equation \[r(s(t))\cos(\beta(t))=-\eta r_{min}\] and since $\beta$ is negative for $t\in \frac{1}{4} \tau(\eta),\frac{3}{4} \tau(\eta)$, we get \[\beta(s^{-1}(z)) = -\arccos\left(\frac{-\eta r_{min}}{r(z)}\right).\]
\begin{align*}
\Rightarrow  &\pi W(\eta)=\int^{\min(s)}_{\max(s)}\frac{\cos\left(-\arccos\left(\frac{-\eta r_{min}}{r(z)}\right)\right)}{r(z)}\frac{1}{\sin\left(-\arccos\left(\frac{-\eta r_{min}}{r(z)}\right)\right)} dz\\
=&\int^{\min(s)}_{\max(s)}\frac{-\eta r_{min}}{r^2(z)}\frac{-1}{\sqrt{1-\left(\frac{-\eta r_{min}}{r(z)}\right)^2}} dz=\int_{\min(s)}^{\max(s)}\frac{1}{r(z)}\frac{1}{\sqrt{\left(\frac{r(z)}{\eta r_{min}}\right)^2-1}} dz
\end{align*}

Because $S_m$ is symmetric with respect to reflecting at the $(x,y)$-plane, we get that

\begin{align*}
&\int_{\min(s)}^{\max(s)}\frac{1}{r(z)}\frac{1}{\sqrt{\left(\frac{r(z)}{\eta r_{min}}\right)^2-1}} dz=2\int_{\min(s)}^{M/4}\frac{1}{r(z)}\frac{1}{\sqrt{\left(\frac{r(z)}{\eta r_{min}}\right)^2-1}}dz\\
=&2\int_{\min(s)}^{\pi/2}\frac{1}{r(z)}\frac{1}{\sqrt{\left(\frac{r(z)}{\eta r_{min}}\right)^2-1}}dz+2\int_{\pi/2}^{M/4}\frac{1}{r(z)}\frac{1}{\sqrt{\left(\frac{r(z)}{\eta r_{min}}\right)^2-1}}dz.
\end{align*}

Because of our assumption that $r(z)=\sin(z)$ for $z\in [0,\pi/2]$, we can explicitly calculate the first integral, and we get
\[\pi W(\eta)=\pi -2\int_{\pi/2}^{M/4}\frac{1}{r(z)}\frac{1}{\sqrt{\left(\frac{r(z)}{\eta r_{min}}\right)^2-1}}dz\]

For each $z \in (\pi/2, M/4)$, the term $\frac{1}{\sqrt{\left(\frac{r(z)}{\eta r_{min}}\right)^2-1}}$ is strictly decreasing in $\eta$ and thus the function $W(\eta)$ is strictly decreasing in $\eta$. By the equation $f(\eta)=LW(\eta)-L$, the same then holds for the map $f$.

\end{proof}

Using this result, we can prove the following theorem about the existence of closed geodesics:

\begin{thm}\label{pq}
For every $(p,q)\in \Z\times \N$ coprime, there exists a unique $S^1$-family of closed geodesics that have Clairaut integral in the range $(-r_{min},r_{min})$ and intersect the Birkhoff annulus $A$ of $M_+$ $q-$times. If $p$ is non-negative, the closed geodesics are $(p+q,q)$-satellites of $E_+$ and $D_+$, if it is negative, the geodesics are $(|p|+q,q)$-satellites of $E_-$ and $D_-$.
\end{thm}

\begin{proof}
We first show $ \lim_{\eta \rightarrow -1}f(\eta)=\infty $. 
Choose a tangent vector $(\Phi(s_{min},\theta),0) $ of $M_+$, where we use parametrization \ref{cut}. The set
\[\Omega:= \{(\Phi(s_{min}+s,\theta),\beta);s,\beta\in (-\epsilon,\epsilon)\}  \]
is a surface transverse to the geodesic flow. Since the geodesic flow is independent of $\theta$, its orbits that intersect $\Omega$ descend to orbits and fixed points of the vectorfield $\sin(\beta)\partial_s +\frac{r^\prime(s)}{r(s)}\cos (\beta)\partial_\beta$ on $\Omega$. From the chosen form of $S$, that is, $r^\prime(s_{min})=0$ and $r^{\prime\prime}(s_{min})>0$, it follows that the origin of this surface is a hyperbolic fixed point of this vector field. The intersection of this surface with the Birkhoff annulus $A$ consists of the points $\{(\Phi(s_{min}+s,\theta),\beta);s=0,\beta\in (0,\epsilon)\}$. Since the origin is a hyperbolic fixed point, geodesics intersecting the Birkhoff annulus at a point $(\Phi(s_{min},\theta),\beta)$ take exponentially longer to leave the neighbourhood \[\{(\Phi(s_{min}+s,\theta),\beta);s,\beta\in (-\epsilon,\epsilon),\theta \in \R/2\pi \Z\}  \]
of $M_+$ as $\beta$ goes to zero. As a consequence, $ \lim_{\eta \rightarrow -1}W(\eta)=\infty $ holds and therefore, by Theorem \ref{tes}, $ \lim_{\eta \rightarrow -1}f(\eta)=\infty $ also holds. The limit $ \lim_{\eta \rightarrow 1}f(\eta)=-\infty $ follows from $f$ being an odd function.\\ From these limits and the fact that $f$ is, by Theorem \ref{dec}, strictly decreasing, it follows that the map $f$ attains every value in $\R$ exactly once.  Therefore, for every $a\in \R$ there exists a unique value $\eta(a)\in (-1,1)$ with $f(\eta(a))=L \cdot a$, where $L$ is the length of $M_+$. Then, given a $(p,q)\in \Z\times \N$ coprime, the point $(x,\eta(\frac{p}{q}))$ is a periodic point of the return map $\rho(x,\eta)=(x+f(\eta),\eta)$ with prime period $q$. Thus, a geodesic starting at $(x,\eta(\frac{p}{q}))$ is closed and intersects $A$ in $q$ different points. The $S^1$-family is then given by varying $x\in \R/L\Z$. If $p$ is greater than zero, $\eta(\frac{p}{q})$ lies in the interval $(-1,0)$, because $f$ is strictly decreasing and $f(0)=0$. Thus, the winding of a closed geodesic $\gamma$ that intersects the annulus at a point $(x,\eta(\frac{p}{q}))$ is, by Theorem \ref{tes} and Definition \ref{windingnumber}, given by $W(\gamma)= q\cdot W(\eta(\frac{p}{q}))=q \cdot (\frac{p}{q}+1)=p+q$. Consequently, by the discussion of Section \ref{linkingnumbers}, the geodesic is a $(p+q,q)$-satellite of $E_+$ and $D_+$. For $p$ less than zero, the lemma follows by the symmetry of $S$.
\end{proof}

\begin{cor}
There exists a subset of the homotopy classes of loops in $S^3\setminus L_1$ with the following properties:
\begin{itemize}
    \item For each homotopy class in this subset, the set of closed Reeb orbits in this homotopy class consists of exactly one $S^1$-family of closed Reeb orbits.
    \item There is a bijective correspondence between the set of coprime $(p,q)\in \Z\times \N$ and this subset.
\end{itemize}
Given a coprime pair $(p,q)\in \Z\times \N$, the corresponding homotopy class of loops in $S^3\setminus L_1$ will be denoted by $y_{(p,q)}$. 
\end{cor}

\begin{proof}
We identify $H_1(S^3\setminus L_1)$ with $\Z^4$ as in Section \ref{lift}. By Theorem \ref{pq} there exists for every $(p,q)\in \Z\times \N$ a unique $S^1$-family $\{\gamma_\theta\}_{\theta\in S^1}$ of geodesics that are $(|p|+q,q)$-satellite of either $E_+$ and $D_+$ or of $E_-$ and $D_-$. Choose a $(p,q)$ and denote by $y_{(p,q)}$ the homotopy class of free loops in $S^3\setminus L_1$ in which the lifts $\tilde{\gamma}_\theta$ of these geodesics to $S^3$ lie. We need to show that there are no other closed Reeb orbits besides the family $\tilde{\gamma}_\theta$ in this homotopy class $y_{(p,q)}$. To do this, we will show that there are no other closed Reeb orbits that represent the same element in $H_1(S^3\setminus L_1)$ as the $\tilde{\gamma}_\theta$ do. According to Lemma \ref{list}, the $\tilde{\gamma}_\theta$ represent the element $(p-q,p-q,-p-q,-p-q)$ or $(\frac{p-q}{2},\frac{p-q}{2},\frac{-p-q}{2},\frac{-p-q}{2})$ if $p>0$ and $(-|p|-q,-|p|-q,|p|-q,|p|-q)$ or $(\frac{-|p|-q}{2},\frac{-|p|-q}{2},\frac{|p|-q}{2},\frac{|p|-q}{2})$ if $p<0$. The list from \ref{list} together with Remark \ref{rema} shows that given a Reeb orbit representing one of these homology classes, its projection onto $S_m$ has to have Clairaut integral in the range $(-r_{min},r_{min})$. This implies that it intersects the Birkhoff annulus at $M_+$ and is therefore part of one of the $S^1$-families. Showing that the lifts of different $S^1$-families represent different elements in $H_1(S^3\setminus L_1)$ will therefore finish the proof. Let $(\tilde{p},\tilde{q}) \in \Z\times \N$, coprime, satisfy $(p,q)\neq (\tilde{p},\tilde{q})$. If $p$ and $\tilde{p}$ have the same sign, the required statement follows from $(|p|-q,-|p|-q)\neq (|\tilde{p}|-\tilde{q},-|\tilde{p}|-\tilde{q})$. If $p$ and $\tilde{p}$ have different signs, the required statement follows from the fact that the first two numbers in $(p-q,p-q,-p-q,-p-q)$ are larger than the last two, while in $(-|p|-q,-|p|-q,|p|-q,|p|-q)$ the first two numbers are smaller than the last two.

\end{proof}

We can provide additional qualitative information about the homotopy classes $y_{(p,q)}$.
\begin{cor}\label{exis}
Given $a,b\in \R, a<b$, then there exists a constant $c$, only depending on $a$ and $b$, such that for every coprime $(p,q)\in \Z\times \N$ with $\frac{p}{q}\in (a,b)$, the closed Reeb orbits in the homotopy class $y_{(p,q)}$ have action bounded by $c \cdot q$. Moreover, $y_{(p,q)}$ is primitive. 
\end{cor}

\begin{proof}
Define $c:=\max\{\tau(\eta(x)); x \in [a, b]\}$. Since a closed geodesic starting at a point $(x,\eta(\frac{p}{q}))$, $x\in \R/L\Z$, intersects the Birkhoff annulus at $M_+$ $q$-times, its length is bounded by $q\cdot c$ and therefore the action of its lift to $S^3$ is bounded by $2\cdot q \cdot c$.
We now show that $y_{(p,q)}$ is a primitive homotopy class. We begin with the case where $p>0$ and $p+q$ is odd. Assume $y_{(p,q)}$ is not primitive. Then the curves in this homotopy class represent the class $(p-q,p-q,-p-q,-p-q)$ in $H_1(S^3\setminus L_1)$, which consequently is also not primitive. Thus, $p-q$ and $-p-q$ have a common divisor different from $\pm 1$, which we denote by $a$. 

\begin{align*}
&a|p-q \: \:\mathrm{and} \:\: a |-p-q\\
\Rightarrow &a|2p \:\:\:\;\;\;\: \mathrm{and}\: \: a|2q
\end{align*}

Since $p$ and $q$ are coprime, this means that $a$ has to be $2$. This gives the desired contradiction, as $p+q$ is odd. \\
The case of $p+q$ even follows similar, as in this case not both $\frac{p+q}{2}$ and $\frac{p-q}{2}$ can be even. The case of $p<0$ follows analogously.

\end{proof}
\begin{defi}
We will denote by $K(a)$ the value of the Clairaut integral of the point $(x,\eta(a))$. By $T^1S_{m(K(a),K(b))}$ we denote the set $K^{-1}(K(a),K(b))$ and by $S^3_{(K(a),K(b))}$ the set $\tilde{K}^{-1}(K(a),K(b))$. Obviously, $S^3_{(K(a),K(b))}$ is the preimage of $T^1S_{(K(a),K(b))}$ under the covering map $\tilde{l}$.

\end{defi}

\section{Homology of the model system}\label{homo}
We would like to compute the cylindrical contact homology of $\lambda_m$, but the contact form $\lambda_m$ does not satisfy the conditions from Section \ref{ch}, namely, it has degenerate Reeb orbits; therefore, its cylindrical contact homology is not well defined. Nevertheless, we will show that it is a Morse-Bott degenerate contact form, allowing us to obtain for any value $S\in \R_{>0}$ a $C^\infty$-close non-degenerate contact form $\lambda_S$, which is non-degenerate up to action $S$.

\subsection{Morse-Bott perturbation}
\begin{defi}
A contact form $\lambda$ on a manifold $M$ is said to be Morse-Bott degenerate if
\begin{enumerate}
\item the action spectrum is discrete,
\item $N_T=\{p\in M \colon \phi_T(p)=p\}$ is a closed, smooth submanifold,
\item the rank of $d\lambda_{|TN_T}$ is locally constant along $N_T$ and
\item the equality $T_p N_T=Ker(d\phi_T(p)-I)$ holds for all $p \in N_T$ .
\end{enumerate}
\end{defi}

\begin{thm}\label{mB}
The restriction of the contact form $\lambda_m$ to $S^3_{(K(a),K(b))}$ is Morse-Bott degenerate for every $a,b\in \R$.
\end{thm}

\begin{proof}
1. We show that for every $T\in \R_{>0}$ there are only finitely many values smaller than $T$ in the action spectrum. For this, choose a number $T>0$.

Given $\eta \in (-1,0)$, the return time to the Birkhoff annulus at $M_+$ is, by Theorem \ref{tes}, given by:
\begin{equation}\label{deriv}
\begin{split}
\tau (\eta)&=F(\eta)-\eta f(\eta)\\
\Rightarrow \tau^\prime(\eta)&=f(\eta)-f(\eta)-\eta f^\prime(\eta)=-\eta f^\prime(\eta)<0
\end{split}
\end{equation}

The return time is therefore, by the fact that $\tau$ is an even function, bounded from below by $\tau(0)=M$ for all $\eta \in (-1,1)$. The lower bound on $\tau(\eta)$ gives an upper bound on the $q\in \N$ satisfying the condition $2q \cdot \tau(\eta)\le T$. The set $Q:=\{q\in \N\colon 2 q \cdot M \le T \}$ is obviously finite. By Theorems \ref{dec} and \ref{pq} every closed orbit in $S^3_{(K(a),K(b))}$ corresponds, up to an $S^1$-action, uniquely to a $(p,q)\in \Z\times \N$ with $\frac{p}{q}\in (a,b)$ and has action $2 q \cdot \tau(\eta(\frac{p}{q}))$. Hence, an orbit can have action less than $T$ only if $q \in Q$.
For each $q \in Q$, the fact that $\frac{p}{q}$ has to be in the interval $(a,b)$ gives a bound on the number of possible $p\in \N$. Therefore, there are only finitely many $(p,q)\in \Z\times \N$ satisfying $q\in Q$ and $\frac{p}{q}\in (a,b)$ and consequently there are only finitely many $S^1$-families of closed periodic orbits with action less than $T$. Since all orbits of one $S^1$-family of closed orbits have the same action, there are only finitely many values less than $T$ in the action spectrum, and thus the action spectrum is discrete. \\
2. By Remark \ref{fibr}, the set $T^1S_{m(K(a), K(b))}$ is foliated by tori given as the level sets of the Clairaut integral $K$. Thus, $S^3_{(K(a), K(b))}$ also gets foliated into tori that arise as the level sets of $\tilde{K}$ and can be written as 
\[S^3_{(K(a),K(b))}=\bigcup_{x\in(a,b)}\tilde{K}^{-1}(K(x)).\]
The level sets $\tilde{K}^{-1}(x)$ are filled by Reeb orbits, which are periodic if and only if $x$ is rational, i.e. if there exist $(p,q)\in\Z\times \N$ coprime with $x=\frac{p}{q}$. These Reeb orbits are exactly the $S^1$-family of closed Reeb orbits in homotopy class $y_{(p,q)}$ and have action $2 q \cdot \tau\left(\eta\left(\frac{p}{q}\right)\right)$. Thus, $N_T$ is the union of all tori $\tilde{K}^{-1}\left(K\left(\frac{p}{q}\right)\right)$, where $p\in \Z, q\in \N$ are coprime and such that $2 q \cdot \tau\left(\eta\left(\frac{p}{q}\right)\right)$ divides $T$. As shown in point $1$, this union is finite.\\
3. The tangent space of $N_T$ is spanned by $\partial_\theta$ and the Reeb vector field $X_{\lambda_m}$.
By formula \ref{ctcform}, $d\lambda_m$ is given by:
\[d\lambda_m=d(r(s)\cos(\beta)d\theta+\sin(\beta)ds)\]

\[\Rightarrow d\lambda_m (\partial_\theta)=-r^{\prime}(s)\cos(\beta)ds+r(s) \sin(\beta)d\beta\neq 0,\]
because the only points $(s,\theta,\beta)$ where both $r^{\prime}(s)\cos(\beta)$ and $r(s) \sin(\beta)$ vanish are the points, where $s$ is a critical point of the function $r$ and $\beta$ is equal to $0$ or $\pi$, but these points are not in $S^3_{(K(a),K(b))}$ as their Clairaut integral is either $\pm r_{min}$ or $\pm r_{max}$. Thus, $\partial_\theta$ is not in the kernel of $d\lambda_m$, but, by definition, $X_{\lambda_m}$ is always in the kernel. Therefore, the rank of $d\lambda_{m_{|TN_T}}$ is equal to one on $TN_T$.\\
4. A basis of the tangent space of $S^3_{(K(a),K(b))}$ is given by the tangentvectors $X_{\lambda_m}, \partial_\theta$ and $\nabla \tilde{K}=r^\prime (s) \cos(\beta)\partial_s-r(s)\sin(\beta)\partial_\beta$, where the first two vectors form a basis of $T_pN_T$. Obviously, $X_{\lambda_m}$ is in $Ker(d\phi_T(p)-I)$. The tangent vector $\partial_\theta$ is also in $Ker(d\phi_T(p)-I)$, as the contact form does not depend on $\theta$. Thus, the inclusion $T_pN_T\subset Ker(d\phi_T(p)-I)$ holds. Let $p$ be a point in $N_T$. It projects to a point in $T^1S_m$ which lies on a geodesic that intersects the Birkhoff annulus at $M_+$ in a point $(x,\eta)$. If $\nabla \tilde{K}$ is in the kernel $Ker(d\phi_T(p)-I)$, then the return time $\tau(\eta)$ has to be constant, which is a contradiction to Equation \ref{deriv}. As a consequence, $\nabla \tilde{K}$ is not in the kernel of $d\phi_T(p)-I$, finishing the proof.

\end{proof}

This allows us to do a Morse-Bott deformation, which we will describe in the proof of the following theorem:

\begin{thm}\label{perturbed}
Given a $(p,q)\in \Z\times \N$ coprime, denote by $N_{(p,q)}$ the unique torus foliated by closed Reeb orbits in homotopy class $y_{(p,q)}$. All of these Reeb orbits have the same action, which will be denoted by $T_{(p,q)}$. By $S_{(p,q)}$ we denote the quotient of $N_{(p,q)}$ under the action induced by the Reeb flow.  \\
Given $a,b\in \R$, for every $S\in \R_{>0}$ and every $\epsilon >0$ there exists a contact form $\lambda_S$ arbitrarily close to $\lambda_m$ in the $C^\infty$-topology such that for all $(p,q)\in \Z\times \N$ coprime with $\frac{p}{q}\in (a,b)$ and $T_{(p,q)}<S$ following conditions hold:

\begin{enumerate}
\item $\lambda_S$ is identical to $\lambda_m$ on a neighbourhood of the link $L_1$.
\item $P^{ y_{(p,q)},\le S}(\lambda_S,L_1)$ consists of exactly two elements which are in a \newline$1-1$-correspondence to the critical points of a perfect Morse function on $S_{(p,q)}$ and have Conley-Zehnder indices differing by one. 
\item The action of the two orbits in $P^{y_{(p,q)},\le S}(\lambda_S,L_1)$ lies in the interval \\$\left(T_{(p,q)}-\epsilon,T_{(p,q)}+\epsilon\right)$.
\item The data $\left(\lambda_S,L_1,y_{(p,q)},S\right)$ satisfies conditions $I)1-4$ and $II)1-2$ from Section \ref{ch}. 

\end{enumerate}

\end{thm}

\begin{proof}

Choose $S \in \R_{>0}$ and $\epsilon > 0$. We first define a function
\[\tilde{g}_S\colon\sqcup_{\{(p,q):T_{p,q}<S\}}N_{(p,q)}\rightarrow \R.\]
As shown in theorem \ref{mB}, there exist only finitely many $(p,q)\in \Z\times\N$ with $T_{(p,q)}<S$ and thus we can extend $\tilde{g}_S$ to a function $g_S\colon S^3 \rightarrow \R$ by cutoff functions. The desired contact form will then be given by $(1+\delta \cdot g_S) \lambda_m$ for a $\delta $ small enough. \\
More detailed, let $(p,q)$ be such that $T_{(p,q)}<S$. The quotient space $S_{(p,q)}$ is diffeomorphic to $S^1$ and after fixing a parametrization 
\begin{equation}\label{paraS}
x\colon \R/2 \pi \Z \rightarrow S_{(p,q)} 
\end{equation}   
we define a function 
\[\overline{g}_{(p,q)}\colon S_{(p,q)}\rightarrow \R\]
\[\overline{g}_{(p,q)}(x)=\cos(x),\]
which is a perfect Morse function.
This function naturally extends to a function $\tilde{g}_{(p,q)}\colon N_{(p,q)} \rightarrow \R$ which is constant on equivalence classes. By performing this construction for every $(p,q)$ with $T_{(p,q)}<S$, we obtain a function $\tilde{g}_S\colon\sqcup_{\{(p,q):T_{p,q}<S\}}N_{(p,q)}\rightarrow \R$ as the collection of these functions.\\
In order to extend this function to $S^3$ by cutoff functions, we need to identify for each such $(p,q)$ the open set $(-\tilde{\epsilon},\tilde{\epsilon}) \times N_{(p,q)}$ with a small neighbourhood $U_{(p,q)}$ of $N_{(p,q)}$ that neither intersects a small neighbourhood of $L_1$ nor any other $U_{(\tilde{p},\tilde{q})}$.

To do this, we first define a basis of $\ker(\lambda_m)$ by:
\begin{equation}\label{basis}
\begin{split}
\xi_1(s,\theta,\beta)&=r(s)\sin(\beta)\cos(\beta) \partial_s+r^{\prime}\cos(\beta)^2 \partial_\beta-\sin^2(\beta) \partial_\theta\\
 \xi_2(s,\theta,\beta) &= \frac{1}{(r^{\prime^2}\cos(\beta)^2+r^2\sin(\beta)^2)}\nabla \tilde{K}(s,\beta)-\frac{\sin r^\prime}{r(r^{\prime^2}\cos(\beta)^2+r^2\sin(\beta)^2)}\partial_\theta\\
 &= \frac{r^\prime(s) \cos(\beta)\partial_s-r(s)\sin(\beta)\partial_\beta}{(r^{\prime^2}\cos(\beta)^2+r^2\sin(\beta)^2)}-\frac{\sin r^\prime}{r(r^{\prime^2}\cos(\beta)^2+r^2\sin(\beta)^2)}\partial_\theta.
 \end{split}
\end{equation}
The vectors are chosen such that $\{\xi_1,\xi_2\}$ forms a basis of $\ker(\lambda_m)$, $\xi_1\in TN_{(p,q)},$ $ d\lambda_m(\xi_1,\xi_2)=1$ and $\frac{\partial \tilde{K}}{\partial \xi_2}=1$.\\
Then we define a flow $\psi\colon(-\tilde{\epsilon},\tilde{\epsilon})\times N_{(p,q)} \rightarrow S^3$ by 
\begin{equation} \frac{\partial \psi}{\partial t}(t,x)=\xi_2(\psi(t,x)).\end{equation}

Since there are only finitely many $(p,q)$ with $\frac{p}{q}\in (a,b)$ and $T_{(p,q)}<S$, there exists a $\tilde{\epsilon}$ small enough such that the flow defines a diffeomorphism $\Psi$ between $(-\tilde{\epsilon},\tilde{\epsilon}) \times N_{(p,q)}$ and a neighbourhood $U_{(p,q)}$ of $N_{(p,q)}$ that neither intersects a small neighbourhood of $L_1$ nor any other $U_{(\tilde{p},\tilde{q})}$. We choose a smooth, compactly supported cutoff function $\beta\colon (-\tilde{\epsilon},\tilde{\epsilon})\rightarrow [0,1]$ with $\beta \equiv 1$ near $0$ and define a function $\hat{g}_{(p,q)}\colon (-\tilde{\epsilon},\tilde{\epsilon}) \times N_{(p,q)}\rightarrow \R$ by
\[\hat{g}_{(p,q)}(t,x)=\beta(t)\tilde{g}_{(p,q)}(x)\]
This allows us to define a function $g_{(p,q)}\colon U_{(p,q)}\rightarrow \R$ by
\[g_{(p,q)}:=\hat{g}_{(p,q)}\circ \Psi^{-1}.\]
As a result, we have $dg_{(p,q)}(\xi_2)(\psi_{t}(x))=\beta^\prime(t)$ for $x\in N_{(p,q)}$.

The desired function $g_S: S^3\rightarrow \R$ is now given by 
\[g_S(x)=\begin{cases}
g_{(p,q)}(x) & \mathrm{if} \: x\in U_{(p,q)} \: \mathrm{for \: some}\: (p,q)\\
0 & \mathrm{else}
\end{cases}\]
Define a function $f_\delta\colon S^3\rightarrow \R$ by \[f_\delta =1+ \delta g_S\] and a $1$-form $\lambda_\delta$ by \[\lambda_\delta =f_\delta \lambda_m.\] The rest of the proof consists of showing that for a $\delta$ small enough, this $1$-form has the properties required by the statement. \\
As a first step, we will compute the Reeb flow of this perturbed contact form. For ease of notation, we will denote the function $f_\delta$ simply by $f$ if the dependency on $\delta$ is not explicitly required. Outside of the union $\sqcup_{\{(p,q):T_{p,q}<S\}}U_{(p,q)}$ the function $f$ is equal to one and therefore the Reeb vector field $X_{\lambda_\delta}$ coincides with the vector field $X_{\lambda_m}$. Now, let $x$ be a point in one of the $U_{(p,q)}$. In the following section, we will suppress the basepoint and assume it to be this point.

Let $X_p$ be a solution of 
\begin{align*}
\lambda_m(X_p)&=0\\
d\lambda_m(X_p)&=\frac{df}{f^2}-\frac{df(X_{\lambda_m})}{f^2}\lambda_m.
\end{align*}

Then, a simple calculation shows that the Reeb vector field of the perturbed contact form $\lambda_\delta$ is given by
\[X_{\lambda_\delta}=\frac{X_{\lambda_m}}{f}+X_p.\]

Using the base of $\ker(\lambda_m)$ defined in \ref{basis} we find that 

\[X_p=\frac{\beta^\prime(t)}{f^2}\xi_1-\frac{df(\xi_1)}{f^2}\xi_2\]

\begin{equation}\label{ver} X_{\lambda_\delta}=\frac{X_{\lambda_m}}{f}+\frac{\beta^\prime(t)}{f^2}\xi_1-\frac{df(\xi_1)}{f^2}\xi_2,
\end{equation}
where $t$ is such that $\psi_t^{-1}(x)\in N_{(p,q)}$.
On $N_{(p,q)}$ we have $df(\xi_1)=-f_\theta$, where $f_\theta$ denotes $\frac{\partial f}{\partial \theta}$. This comes from the fact that on $N_{(p,q)}$ it holds that $df(X_{\lambda_m})=0$ and $\xi_1=r(s)\cos(\beta)X_{\lambda_m}-\partial_\theta$.

We can now show that the conditions of the theorem are satisfied for a $\delta$ small enough.
\begin{enumerate}
\item By construction $\lambda_S$ is identical to $\lambda_m$ near $L_1$, since $f_\delta $ is equal to $1$ near $L_1$.
\item 
\begin{enumerate}
\item \label{2a} At the two Reeb orbits on $N_{(p,q)}$ that correspond to the critical values of the perfect Morse function $\overline{g}_{(p,q)}$ the partial derivative $f_\theta$ vanishes and therefore the Reeb vector field $X_{\lambda_\delta}$ is a multiple of the Reeb vector field $X_{\lambda_M}$. As a consequence, both orbits remain Reeb orbits under the perturbation, now with action $(1+\delta)T_{(p,q)}$ and $(1-\delta)T_{(p,q)}$, respectively. Since the inequality $T_{(p,q)}<S$ holds by assumption, the inequality $(1+\delta)T_{(p,q)}< S$ holds for $\delta$ small enough. The orbit corresponding to the maximum of the perfect Morse function $\overline{g}_{(p,q)}$ will be denoted by $P_{max}$, the one corresponding to the minimum will be denoted by $P_{min}$.
\item \label{2b} Let $p$ be a point in the image of $P_{max}$ or $P_{min}$. In the following calculation, we will assume the base point to be $p$. Using the global trivialization \ref{basis} of $\xi$ and denoting by $\phi_t$ the Reeb flow of $X_{\lambda_\delta}$, write $d\phi_t$ as a matrix to get the equation 
\[\frac{\partial}{\partial t} d\phi_t=DX_{\lambda_\delta}d\phi_t=\left(D\left(\frac{X_{\lambda_m}}{f}\right)+DX_p\right)d\phi_t.\]

with
\[D\left(\frac{X_{\lambda_m}}{f}\right)=\begin{pmatrix}
0& c(p,q)\\
0&0\\
\end{pmatrix},
DX_{p}=\begin{pmatrix}
0&0\\
\frac{\partial} {\partial \xi_1} \frac{f_{\theta}}{f^2}&0\\
\end{pmatrix},\]

with $c(p,q)\in\R\neq 0$. Since the basepoint is in the image of $P_{max}$ or $P_{min}$ it holds that $\frac{\partial} {\partial \xi_1} \frac{f_{\theta}}{f^2}=\frac{\partial }{\partial \theta}\frac{f_{\theta}}{f^2}= \frac{\pm \delta }{(1\mp \delta)^2} $. Independent of the sign of $c$, these two matrices show that the difference between the Conley-Zehnder indices of $P_{max}$ and $P_{min}$ is always equal to one, see for example \cite{GCZ}.

Additionally, both orbits are non-degenerate.
\item All other periodic Reeb orbits of $X_{\lambda_m}$ on $N_{(p,q)}$ are not periodic orbits of $X_{\lambda_\delta}$, because the $\xi_2$-component of $X_{\lambda_m}$ and therefore also the $\nabla \tilde{K}$-component of $X_{\lambda_m}$ is a constant different from zero along the other orbits. 
\item \label{2d} We have already shown that $P^{\le S, y_{(p,q)}}(\lambda_S,L_1)$ contains at least two elements; we now show that it does not contain any more elements. We prove this by showing that for $\delta$ small enough, the contact form $\lambda_\delta$ has no new closed Reeb orbits with action less than $S$ in each $U_{(p,q)}$. Assume the contrary: that there exists a sequence $\delta_i$ converging to $0$ such that for every $\delta_i$ there exists a closed $\lambda_{\delta_i}$-Reeb orbit $\gamma_{\delta_i}$ different from $P_{max}$ and $P_{min}$ with action less than $S$ in $U_{(p,q)}$.
Then the sequence $\gamma_{\delta_i}$ has a convergent subsequence. 
The limit $\gamma$ of this subsequence has to be a closed orbit on $N_{(p,q)}$, as it is a closed $X_{\lambda_m}$-Reeb orbit and the only closed $X_{\lambda_m}$-Reeb orbits in $U_{(p,q)}$ with action less than $S$ are subsets of $N_{(p,q)}$. Furthermore, it has to be $P_{max}$ or $P_{min}$, since all other orbits have a small neighbourhood in which the $\xi_2-$ and thus the $\nabla \tilde{K}$-component of the Reeb vector field is different from zero for all $\delta$. Using the normal bundle of $P_{max}$ given by $\{\xi_1,\xi_2\}$ identify a small neighbourhood of $P_{max}$ with $S^1\times (-\tilde{\epsilon},\tilde{\epsilon})\times (-\tilde{\epsilon},\tilde{\epsilon})$. Since the Reeb flow $X_{\lambda_\delta}$ is independent of the $S^1$-direction it projects to a flow on $(-\tilde{\epsilon},\tilde{\epsilon})\times (-\tilde{\epsilon},\tilde{\epsilon})$ and orbits of $X_{\lambda_\delta}$ project to orbits and fixed points of this flow. The origin, corresponding to $P_{max}$, is a fixed point of this flow. Furthermore, using the characteristic equation of $DX_{\lambda_\delta}$, it is possible to see that one of $P_{max}$ and $P_{min}$ is hyperbolic, while the other one is elliptic; therefore, the origin has to be either a hyperbolic or an elliptic fixed point. If it is a hyperbolic fixed point, the projection of the family $\gamma_{\delta_i}$ cannot converge to it, as otherwise it would bound a disc containing either no singular point or only a hyperbolic singularity, both contradictions to $\chi(D)=1$ by the Poincaré-Hopf theorem. If the fixed point is elliptic, the projection of the family $\gamma_{\delta_i}$ cannot converge to it, because the eigenvalues of $DX_{\lambda_\delta}$ are proportional to $\sqrt{\delta}$ and therefore the projection of the $\gamma_{\delta_i}$ would take a long time to close up, contradicting the fact that the periods are bounded uniformly by $S$. For this reason, the family $\gamma_{\delta_i}$ cannot converge to $P_{max}$, and the same argument for $P_{min}$ shows that the family can also not converge to $P_{min}$. This is the required contradiction.

\end{enumerate}

\item As shown in point \ref{2a} the action of the two orbits in $P^{\le S, y_{(p,q)}}(\lambda_S,L_1)$ is equal to $(1+\delta)T_{(p,q)}$ and $(1-\delta)T_{(p,q)}$. For $\delta$ small enough these values obviously lie in the interval $(T_{(p,q)}-\epsilon,T_{(p,q)}+\epsilon)$.
\item 
\begin{enumerate}
\item I) 1 holds obviously. 
\item I) 2 holds because of Corollary \ref{exis}.
\item I) 3 holds because of the definition of the homotopy class $y_{(p,q)}$, the linking numbers from \ref{lkn} and \ref{linkingnumbersorbits} and the fact that $p$ and $q$ are coprime. 
\item I) 4 holds because of the linking numbers from \ref{linkingnumbersorbits}.
\item II) 1: For $\delta$ small enough there exists no closed Reeb orbit $\gamma$ in $S^3\setminus L_1$ with action less than $S$ that is contractible in $S^3\setminus L_1$. This follows from point \ref{2d} and the fact that $\lambda_m$ has no closed orbits contractible in $S^3\setminus L_1$.
\item II) 2 holds because of \ref{2b}.

\end{enumerate}

\end{enumerate}

Thus, there exists a $\delta_0$ small enough such that the contact form $\lambda_S:=\lambda_{\delta_0}$ satisfies the conditions of the theorem.

\end{proof}

\subsection{Cylindrical contact homology of the model system}
Theorem \ref{perturbed} tells us that there are exactly two orbits in $P^{y_{(p,q)},S}(\lambda_S,L_1)$. The final step to determine the homology $HC_\ast^{y_{(p,q)}\le S}(\lambda_S, L_1)$ is to define a regular almost complex structure $J$ and to calculate the differential map, which amounts to determining the moduli space $M_{J}^{y_{(p,q)},\le S}(P_{max}, P_{min})$.

\begin{cor}\label{results}
Let $(a,b)\subset \R$ be an interval and $S\in \R_{>0}$. Then there exists a regular almost complex structure $J$ on $\R\times S^3$ such that for all $(p,q)\in \Z\times \N$ coprime with $\frac{p}{q}\in (a,b)$ and $T_{(p,q)}<S$, the moduli space $M_{J}^{y_{(p,q)},\le S}(P_{max},P_{min};L_1)$ consists of exactly two elements. As a result, \[HC_\ast^{y_{(p,q)}\le S}(\lambda_S,L_1)=H_{\ast -c}(S^1,\Z/2\Z)\] holds for some $c\in \Z$.
\end{cor}

\begin{proof}

We define an almost complex structure $J$ on $\R\times S^3$ by
\begin{align*}
J(\xi_1) &= \xi_2 \\
J(\partial_t) &=X_{\lambda_S},
\end{align*}
where $\partial_t$ denotes the $\R$-direction. A simple calculation shows that this almost complex structure is $d(e^t\lambda_S)$-compatible.\\
Next, we explicitly define two $J$-holomorphic cylinders. In order to do this, let $\nabla_Jf\in \Gamma (TS^3)$ be a solution to 
\[ \lambda_S(\nabla_Jf)\lambda_S+ d\lambda_S(...,J \nabla_Jf)=df.\]
A simple calculation shows that $JfX_p=\nabla_Jf$ holds. Thus, the vector field $\nabla_Jf$ is equal to $-\frac{f_\theta}{f} \xi_1$ on $N_{(p,q)}$. Using the parametrization \ref{paraS} of $S_{(p,q)}$ and the fact that $\xi_1=r(s)\cos(\beta)X_{\lambda_m}-\partial_\theta$ holds, we get that  $\nabla_Jf$ descends to the vector field $-\frac{\sin(x)}{1+\delta_0 \cos(x) }\partial_x$ on $S_{(p,q)}$.
 
Let $\eta(s)$ be a flowline of $\nabla_Jf$, i.e. $\eta^\prime(s)=\nabla_Jf(\eta(s))$, such that $\eta(0)$ is a point $p$ on $N_{(p,q)}$ that descends to the point $x=\frac{\pi}{2}$ in $S_{(p,q)}$. Then it holds that $\lim_{s\rightarrow \infty}\eta(s) \in \mathrm{Im}(P_{max})$ and $\lim_{s\rightarrow -\infty}\eta(s) \in \mathrm{Im}(P_{min})$. Therefore, the following $J$-holomorphic cylinder $u_1$ is asymptotic to $P_{max}$ at $+\infty\times S^1$ and to $P_{min}$ at $-\infty\times S^1$. Define
\[u_1\colon \R \times \R/\Z\rightarrow \R\times S^3\]
\[ u_1(s,t)=(a(s), \phi_{T\cdot t}(\eta(T\cdot s))),\]
where $\phi$ denotes the Reeb flow of $\lambda_m$ and the function $a$ is defined by
\[a(s)=T\int_0^s f(\eta(T\cdot \tau))d\tau.\]

We will denote by $\tilde{u}_1$ the $S^3$-component of $u_1$, i.e. $\tilde{u}_1=\phi_{T\cdot t}(\eta(T\cdot s))$.

We will now show that this map is actually a $J$-holomorphic map. For ease of notation, we will denote in the following calculation $\phi_{T\cdot t}(\eta(T\cdot s))$ simply by $\tilde{\phi}$. Then we get the partial derivatives
\begin{align*}
\partial_su_1(s,t)&=(Tf(\eta(T\cdot s)),T\cdot D\phi_{T\cdot t}\nabla_Jf(\eta(T\cdot s)))\\
&=(Tf(\eta(T\cdot s)),T\cdot \nabla_Jf(\tilde{\phi}))\\
&=(Tf(\eta(T\cdot s)),-T\cdot \frac{f_\theta(\tilde{\phi})}{f(\tilde{\phi})}\xi_1(\tilde{\phi})),
\end{align*}
where the second equality holds, because $T_p N_T=\mathrm{ker}(d\phi_T(p)-I)$ holds for all $p \in N_T$.
\begin{align*}
\partial_tu_1(s,t)&=(0,TX_{\lambda_m}(\tilde{\phi}))=(0,Tf(\tilde{\phi})X_{\lambda_S}(\tilde{\phi})-Tf(\tilde{\phi})X_p(\tilde{\phi}))\\
&=(0,Tf(\tilde{\phi})X_{\lambda_S}(\tilde{\phi})-T\frac{f_\theta(\tilde{\phi})}{f(\tilde{\phi})}\xi_2(\tilde{\phi}))\\
\Rightarrow J\partial_tu_1(s,t)&=(-Tf(\tilde{\phi}),T\cdot\frac{f_\theta(\tilde{\phi})}{f(\tilde{\phi})}\xi_1(\tilde{\phi}))\\
&=(-Tf(\eta(T\cdot s)),T\cdot \frac{f_\theta}{f}(\tilde{\phi})\xi_1(\tilde{\phi})),
\end{align*}

where the last equality holds, because on $N_{(p,q)}$ the function $f$ is invariant under the Reeb flow of $X_{\lambda_m}$. This shows that $\partial_su_1(s,t)+J\partial_tu_1(s,t)=0$ and therefore $u_1$ is $J-$holomorphic.

We define a second $J$-holomorphic cylinder $u_2$ by choosing a flowline $\eta_2(s)$ such that $\eta(0)$ is a point $p$ on $N_{(p,q)}$ that descends to the point $x=\frac{3\pi}{2}$ in $S_{(p,q)}$ and then doing the same construction as for $u_1$.

Similar to \cite{UHAM} section 4, it is possible to use Siefring's intersection theory (\cite{RSI}
\cite{MR2456182}) to rule out the existence of a third element in $M_{J}^{y_{(p,q)},\le S}(P_{max}, P_{min}, L_1)$.

The last step we need to prove is that the two holomorphic curves are cut out transversely.\\
For this, we use an automatic transversality result by Wendl:
\begin{thm}[Wendl \cite{C}.]
If $ind(u)>-2+2g + \#\Gamma_0+\#\pi_0(\partial \Sigma)+2 Z(du)$, then $u$ is regular.
\end{thm}

We compute the different values in the expression:
\[ind(u)=\mu_{CZ}(P_{max})-\mu_{CZ}(P_{min})=1\] 
From this equation it also follows that precisely one of the orbits is even and therefore $\#\Gamma_0$ is equal to one. The domain of the $u_i$ is $S^2$ minus two points; thus, the genus $g$ is zero and there is no boundary. Finally, the differential of $u_i$ never vanishes, and $Z(du)$ is therefore equal to zero. Plugging these values into the formula results in:
\[1>-2 +0 + 1 +0 +0\]
This is obviously true, which proves that $u_1$ and $u_2$ are regular.

\end{proof}

\section{Main theorem}
Now that we have all the necessary prerequisites, we will prove the main theorem in this section.
\begin{thm}\label{main}
Let $g$ be a reversible Finsler metric on $S^2$ and define $P^t(g)$ by \[P^t(g):=\# \{\gamma:\gamma \mathrm{\:is\: a \:closed \: prime \: geodesic  \: with \: length \: less\: than \:}t\}.\]
Then 
\[\liminf_{t\rightarrow \infty}\frac{\log(P^t(g))}{\log(t)}\ge 2\]
holds.

\end{thm}
The first step will be to prove the following theorem.
\begin{thm}\label{central}
Let $g$ be a reversible Finsler metric on $S^2$ that has two disjoint, simple closed geodesics.
Then 
\[\liminf_{t\rightarrow \infty}\frac{\log(P^t(g))}{\log(t)}\ge 2\]
holds.
\end{thm}

Given a reversible Finsler metric $g$ on $S^2$, we denote by $\lambda$ the lift of the Hilbert contact form associated to $g$ to $S^3$. The two disjoint, simple closed geodesics together with the two curves created by reversing the orientation form a link $K_0$ that is isotopic to the link $K_1$ from Section \ref{lift}. Thus, the lift of $K_0$ is transversely isotopic to the link $L_1$ of the model system $\lambda_m$. Then Theorem \ref{central} follows from Theorem \ref{number} and the following theorem.

\begin{thm}\label{degenprim}
Let $\lambda$ be a tight contact form on $S^3$ that realizes a link $L$ that is transversely isotopic to $L_1$ as a link of closed Reeb orbits. Let $(a,b)\subset \R$ be an interval. Then there exists a $c\in \R_{>0}$, only depending on the choice of $(a,b)$, such that for each $(p,q)\in \Z\times \N$ coprime satisfying $\frac{p}{q}\in (a,b)$ there exists a closed Reeb orbit in homotopy class $y_{(p,q)}$ whose action is bounded from above by $q\cdot c$.
\end{thm}

Following Hryniewicz-Momin-Salomão \cite{UHAM}, we will use a neck stretching argument to prove Theorem \ref{degenprim}.  More precisely, let $(a,b)\subset \R$ be an interval and $S\in \R_{>0}$ a value. Furthermore, let $\lambda_S$ be the model contact form on $S^3$ from Theorem \ref{perturbed}. Since $S^3$ is compact and both $\lambda$ and $\lambda_S$ are tight, we can choose values $c_\pm \in \R_{>0}$ such that $c_-\lambda_S \prec \lambda \prec c_+\lambda_S$ holds. Additionally, we choose almost complex structures $J_+\in \mathcal{J}_{reg}(c_+\lambda_S)$, which induces an almost complex structure $J_-\in \mathcal{J}_{reg}(c_-\lambda_S)$ as in Section \ref{comm}, $J\in \mathcal{J}_{reg}(\lambda)$, $J_1 \in \mathcal{J}_{reg}(J_-,J;L_1)$ and $J_2 \in \mathcal{J}_{reg}(J, J_+;L_1)$. For each $R\in \R_{>0}$ the almost complex structures $J_1$ and $J_2$ define another almost complex structure $J_R$ as in Section \ref{split}. We can now prove the existence of $J_R$-holomorphic cylinders, which we will subsequently use in the neck stretching process.

\begin{thm}\label{exist}
Let $(a,b)\subset \R$ be an interval and $S\in \R_{>0}$ a value. Then for all $(p,q)\in \Z\times \N$ coprime, with $\frac{p}{q}\in (a,b)$ and $T_{(p,q)}<S$, and every $R\in \R_{> 0}$ there exists a $J_R$-holomorphic cylinder which does not intersect $\tau^{-1}(L_1)$ and is asymptotic to a Reeb orbit $\gamma_+ \in P^{y_{(p,q)},\le S}(c_+ \lambda_S;L_1)$ at its positive puncture and to a Reeb orbit $\gamma_- \in P^{y_{(p,q)},\le S}(c_- \lambda_S;L_1)$ at its negative puncture.
\end{thm}
\begin{proof}
Choose a $(p,q)\in \Z\times \N$ coprime with $\frac{p}{q}\in (a,b)$ and $T_{(p,q)}<S$. Assume the contrary, that there exists an $R\in \R$ for which there exists no such $J_R$-holomorphic cylinder. Define $c=\frac{c_-}{c_+}$. As shown in Section \ref{homo}, the homologies $HC_\ast^{y_{(p,q)},\le S}(c_\pm\lambda_S;L_1)$ and $HC_\ast^{y_{(p,q)},\le S/c}(c_+\lambda_S;L_1)$ are well defined and up to a common shift in the degree equal to $H_{\ast}(S^1,\Z/2\Z)$. Furthermore, as shown in Section \ref{comm}, they are generated by the same orbits and the differential counts the same holomorphic cylinders. Thus, the inclusion map $$i_\ast\colon C_\ast^{y_{(p,q)},\le S}(c_+\lambda_S;L_1)\rightarrow C_\ast^{y_{(p,q)},\le S/c}(c_+\lambda_S;L_1)$$ is nontrivial and the map $$j_\ast \colon C_\ast^{y_{(p,q)},\le S/c}(c_+\lambda_S;L_1)\rightarrow C_\ast^{y_{(p,q)},\le S}(c_-\lambda_S;L_1)$$ is an isomorphism on the level of homology. Consequently, the map \[j_\ast \circ i_\ast\colon C_\ast^{y_{(p,q)},\le S}(c_+\lambda_S;L_1)\rightarrow C_\ast^{y_{(p,q)},\le S}(c_-\lambda_S;L_1)\] is nontrivial. By Theorem \ref{gleich}, this map is chain homotopic to the map $\Phi(J)$ from Section \ref{chain map} for any $J\in \mathcal{J}_{reg}(J_-,J_+:L_1)$ and thus these maps are also nontrivial.

As shown in Section \ref{split}, $J_R$ is biholomorphic to an almost complex structure $\tilde{J}_R\in \mathcal{J}(J_-,J_+:L_1)$. By our assumption, there exists no $J_R$-holomorphic cylinder and therefore likewise no $\tilde{J}_R$-holomorphic cylinder. As a consequence, $\tilde{J}_R$ is automatically in $\mathcal{J}_{reg}(J_-,J_+:L_1)$ and furthermore $\Phi(\tilde{J}_R)$ is trivial. This gives the required contradiction.
\end{proof}

\begin{proof}[Proof of theorem \ref{degenprim}]
Let $(p,q)\in \Z\times N$ be coprime with $\frac{p}{q}\in (a,b)$. Choose a value $S>T_{(p,q)}$. A neck stretching argument similar to \cite{UHAM} for the family of cylinders from Theorem \ref{exist} shows the existence of an element of $P^{y_{(p,q)},\le S}(\lambda,L_1)$. Since $T_{(p,q)}$ is uniformly bounded in $q$ by Corollary \ref{exis}, this completes the proof.
\end{proof}

\begin{proof}[Proof of theorem \ref{central}]
Choose $a,b\in\R$ and denote by $\lambda$ the lift of the Hilbert contact form associated to $g$ to $S^3$. The contact form $\lambda$ is tight and, as mentioned previously, it realizes a link $L$ that is transversely isotopic to $L_1$ as a link of Reeb orbits. Therefore, we can apply Theorem \ref{degenprim}. Thus, for every $(p,q)\in \Z\times \N$ coprime with $\frac{p}{q}\in(a,b)$ there exists a periodic $\lambda$-Reeb orbit whose action is uniformly bounded in $q$. The projection of each such periodic Reeb orbit to $S^2$ is a closed geodesic whose minimal period is less than or equal to the action of the Reeb orbit. Thus, for every such $(p,q)$, there also exists a closed geodesic with length uniformly bounded in $q$. Then, Theorem \ref{number} finishes the proof.
\end{proof}

\begin{proof}[Proof of theorem \ref{main}]
If the metric $g$ is Riemannian, the works of Lyusternik- \newline Schnirelmann \cite{MR0029532},\cite{lyusternik1929probleme} and Bangert \cite{B1} show that the metric satisfies one of the following two cases:
\begin{enumerate}
\item There exist two disjoint, closed, simple geodesics.
\item There exists a simple closed geodesic whose Birkhoff annulus is a Birkhoff section and gets intersected by at least one different closed geodesic.
\end{enumerate}
The work of De Philippis, Marini, Mazzucchelli and Suhr \cite{MR4696545} and Contreras, Knieper, Mazzucchelli and Schulz \cite{CKMSS} shows that this case distinction also holds for reversible Finsler metrics. Thus, we can always assume that our given metric $g$ satisfies one of the two cases.
In the first case, the theorem was proven in \ref{central}.\\
In the second case, there exists a well-defined return map $r\colon A \rightarrow A$, which has a periodic point, and Theorem \ref{FI} finishes the proof.
\end{proof}

\printbibliography

\end{document}